\PassOptionsToPackage{numbers,sort}{natbib}
\documentclass[graybox,natbib,envcountsame]{svmult}

\usepackage{mathptmx}
\usepackage{helvet}
\usepackage{courier}
\usepackage{makeidx}
\usepackage{graphicx}
\usepackage{multicol}
\usepackage[bottom]{footmisc}

\usepackage[utf8]{inputenc}
\usepackage{amssymb}
\usepackage{amsmath}   
\usepackage[
  bookmarks=false, 
 hidelinks,
  ]{hyperref}
\usepackage{newtxtext,newtxmath} 
\usepackage[mathscr]{eucal}
\usepackage{xspace}
\usepackage{mathtools}
\usepackage{paralist}
\usepackage{float}
\usepackage{tikz}
  \usetikzlibrary{graphs,graphs.standard}

\newenvironment{discussedreferences}%
{\begin{list}{}{%
      \setlength{\itemindent}{0mm}%
      \setlength{\listparindent}{0mm}%
      \setlength{\topsep}{1\baselineskip}%
      \setlength{\rightmargin}{0mm}%
      \setlength{\leftmargin}{.1\linewidth}}%
      \setlength{\itemsep}{0pt}%
      \setlength{\topsep}{0pt}%
      \setlength{\parsep}{0pt}%
      \setlength{\parskip}{1ex plus 1ex}%
\item[]\small\raggedright}
{\unskip\end{list}}

\renewcommand{\descriptionlabel}[1]{\hspace{\labelsep}\textit{#1:}}

\newcommand{\arXiv}[1]{\href{http://arxiv.org/abs/#1}{arXiv:#1}}%

\newcommand{\cP}{\ensuremath{\mathscr{P}}\xspace}
\newcommand{\cNP}{\ensuremath{\mathscr{NP}}\xspace}
\let\emptyset\varnothing
\let\phi\varphi

\DeclareMathOperator{\inj}{inj}
\DeclareMathOperator{\vol}{vol}
\DeclareMathOperator{\conv}{conv}
\DeclareMathOperator{\tr}{tr}
\DeclareMathOperator{\Name}{Name}
\DeclareMathOperator{\TSP}{TSP}
\DeclareMathOperator{\TSPLP}{TSPLP}
\DeclareMathOperator{\TSGLP}{TSGLP}
\DeclareMathOperator{\STAB}{STAB}
\DeclareMathOperator{\QSTAB}{QSTAB}
\DeclareMathOperator{\QSTAG}{QSTAG}
\let\TH\relax
\DeclareMathOperator{\TH}{TH}

\spnewtheorem*{nntheorem}{Theorem.}{\bfseries}{\itshape}
\spnewtheorem*{nncorollary}{Corollary.}{\bfseries}{\itshape}
\spnewtheorem*{nnlemma}{Lemma.}{\bfseries}{\itshape}
\spnewtheorem*{noheadingthm}{\unskip}{\bfseries}{\itshape}

\spnewtheorem*{SOPT}{The Strong Optimization Problem (SOPT).}{\bfseries}{\itshape}
\spnewtheorem*{SSEP}{The Strong Separation Problem (SSEP).}{\bfseries}{\itshape}
\spnewtheorem*{SMEM}{The Strong Membership Problem (SMEM).}{\bfseries}{\itshape}
\spnewtheorem*{WOPT}{The Weak Optimization Problem (WOPT).}{\bfseries}{\itshape}
\spnewtheorem*{WSEP}{The Weak Separation Problem (WSEP).}{\bfseries}{\itshape}
\spnewtheorem*{WMEM}{The Weak Membership Problem (WMEM).}{\bfseries}{\itshape}
\spnewtheorem*{SFM}{Submodular Function Minimization.}{\bfseries}{\itshape}

\graphicspath{{./0-img/}}

\begin{document}

\title*{The Mathematics of László Lovász}
\author{Martin Grötschel \and Jaroslav Nešetřil}
\institute{Martin Grötschel \at
  Technische Universität Berlin,
  Mathematisches Institut, 
  Straße des 17.~Juni~135, 
  10623~Berlin,
  Germany,
  \email{groetschel@bbaw.de} \and %
  Jaroslav Nešetřil \at 
  Computer Science Institute of Charles University,
  Faculty of Mathematics of Physics, 
  Charles Univ\-ersity,
  Malostransk\' e n\' am.~25, 
  11800~Praha~1, 
  Czech Republic,
  \email{nesetril@iuuk.mff.cuni.cz}}
\maketitle

\setcounter{minitocdepth}{2}
\dominitoc

\vspace{1\baselineskip}

\abstract{This is an exposition of the contributions of László Lovász to
  mathematics and computer science written on the occasion of the bestowal of
  the Abel Prize~2021 to him. Our survey, of course, cannot be exhaustive. We
  sketch remarkable results that solved well-known open and important problems
  and that -- in addition -- had lasting impact on the development of
  subsequent research and even started whole new theories. Although discrete
  mathematics is what one can call the Lovász home turf, his interests were,
  from the beginning of his academic career, much broader. He employed algebra,
  geometry, topology, analysis, stochastics, statistical physics, optimization,
  and complexity theory, to name a few, to contribute significantly to the
  explosive growth of combinatorics; but he also exported combinatorial
  techniques to many other fields, and thus built enduring bridges between
  several branches of mathematics and computer science. Topics such as
  computational convexity or topological combinatorics, for example, would not
  exist without his fundamental results. We also briefly mention his
  substantial influence on various developments in applied mathematics such as
the optimization of real-world applications and cryptography.}

\section{Introduction}
\label{sec:GN1}

László Lovász was born in 1948 in Budapest. Laci, as he is called by his
friends, attended the Fazekas Mihály Gimnázium in Budapest, a special school
for mathematically gifted students and a fertile ground of world-class
mathematicians. Katalin Vesztergombi, his wife since 1969, was one of his
classmates. Laci's outstanding talent became visible at very young age. He won,
for example, several mathematics competitions in Hungary and also won three
gold medals in the International Mathematical Olympiad.

\begin{figure}[b]
\centering
\includegraphics[width=.7\linewidth]{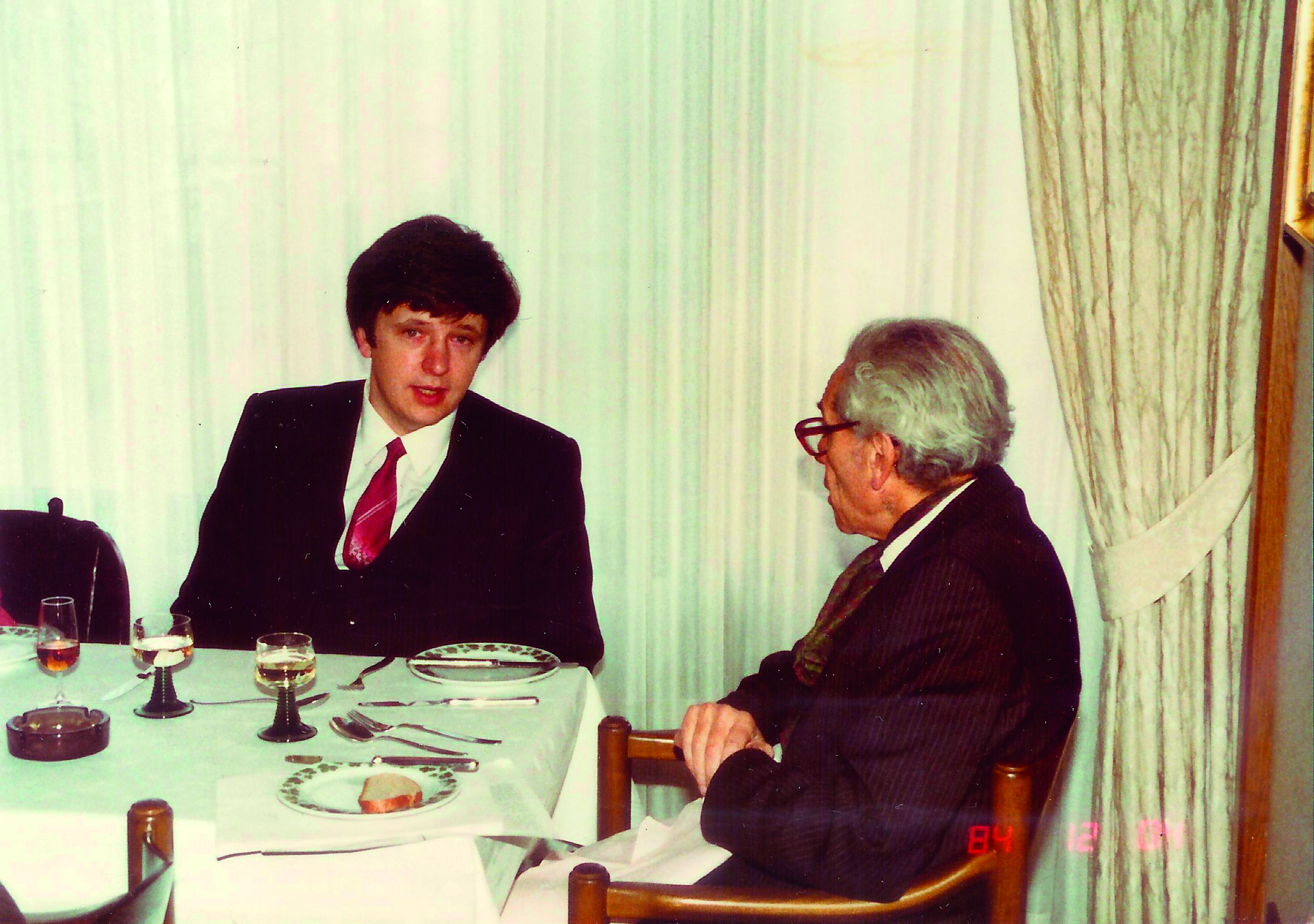}

\caption{Lovász and Erdős at dinner in 1977 (Photo: Private)}
\label{fig:GN1}

\end{figure}

Lovász studied mathematics at Eötvös Loránd University (ELTE). He received \mbox{--
with} Tibor Gallai as his mentor -- his first doctorate (Dr.~Rer.~Nat.) degree
from ELTE in 1971, the Candidate of Sciences (C.\,Sc.) degree in 1970 and his
second doctorate (Dr.\,Math.\,Sci.) degree in 1977 from the Hungarian Academy
of Sciences. Of great influence for his scientific growth was the outstanding
Hungarian combinatorial school (e.g., T.~Gallai, A.~Hajnal, A.~R\'enyi, M.~Simonovits,
V.\,T.~Sós, P.~Turán, and foremost P.~Erdős).

In 1971 Lovász started his professional career as a research associate at ELTE.
From 1975 to 1982 he was Docent, later Professor and Chair of Geometry at
József Attila University, Szeged; 1983--1993 Chair of Computer Science at ELTE;
1993--1999 Professor of Computer Science at Yale University; and 1999--2006
Senior Researcher, Microsoft Research, Redmond. In~2006 Lovász returned to his
hometown Budapest as a Professor and Director of the Mathematical Institute at
ELTE from which he retired in~2018. In 2020 he joined the Alfréd Rényi
Institute of Mathematics.  Lovász served the International Mathematical Union
as its President from~2007 to~2010 and the Hungarian Academy of Sciences as its
President from~2014 to~2020 during demanding times.

\begin{figure}[t]
\centering
\includegraphics[width=.7\linewidth]{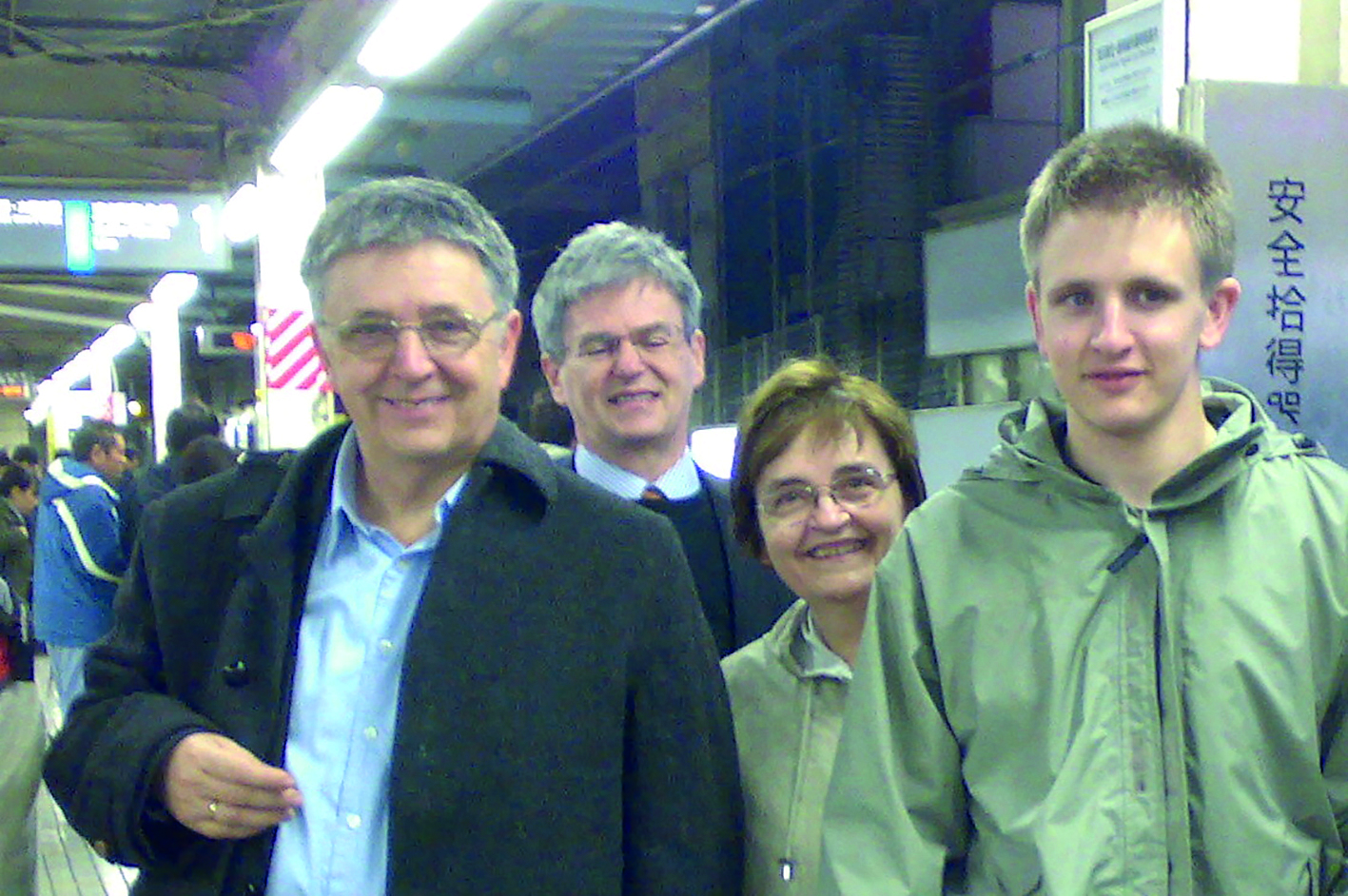}

\caption{In a Tokyo subway station on the way to the Kyoto Prize ceremony:
Laci, Kati, and son Laci M. Lovász, Andr\'as Frank in the back (Photo: Private)}
\label{fig:GN2}
\end{figure}

Among the institutions Lovász visited for extended periods of time are
Vanderbilt University, University of Waterloo, Universität Bonn, University of
Chicago, Cornell University, Mathematical Sciences Research Institute in
Berkeley, Princeton University, Princeton Institute for Advanced Study, and ETH
Zürich. Five universities bestowed special professorships upon him, he received
six honorary degrees and countless high-ranking honors and distinctions,
including the Kyoto Prize~2010, see~Fig.~\ref{fig:GN2}.

Like every scientific discipline, mathematics has become a field with a
large number of specializations. The Mathematics Subject Classification
(MSC~2020) with its 63~first-level areas and 6,006~specific research
areas is a witness of this development. Today, no mathematician has a
full understanding of all the mathematical branches. But there are still
a few people with broad mathematical knowledge, deep command of their
fields of special interest, and the ability to build bridges by
transferring results and techniques between fields to expand the
mathematical toolboxes and open up new research areas. One of these rare
persons is László Lovász. In fact, quite fittingly, two volumes
published in his honor at special occasions were entitled \emph{Building
Bridges}, see~\cite{GroetschelKatona2008} and~\cite{BaranyKatonaSali2019}.

Laci's mathematical roots are in combinatorics. But he vastly expanded his
reach by employing combinatorial methods in other mathematical fields and
bringing, in return, tools from geometry, topology, algebra, analysis,
probability theory, information theory, optimization, and even ideas from
physics into combinatorics. His deep interest in algorithms led to major
advances in modern complexity theory. In his work, Lovász established profound
connections between discrete mathematics and computer science. This is
reflected in the statement that the Norwegian Academy of Science and Letters
issued in its announcement of the award of the Abel Prize 2021 to him and Avi
Wigderson
\begin{quote}
for their foundational contributions to
theoretical computer science and discrete mathematics, and their leading role
in shaping them into central fields of modern mathematics.
\end{quote}
At the end of the 1960s and the beginning of the 1970s, graph theory, discrete
mathematics, combinatorics, and theoretical computer science were considered
peripheral fields of mathematics. This changed completely during Lovász's
lifetime. They became central parts of modern mathematics for many reasons. The
tremendous development of computer technologies is the most obvious one.
Essential factors were also the high quality of the research and the results in
these areas and their wide applicability. The solutions of many problems
arising in industry, society, other sciences, even in other fields within
mathematics critically depend on theories and algorithms invented in discrete
mathematics. Many mathematicians and computer scientists contributed to this.
László Lovász undoubtedly was and still is one of the key players in this
development.

There are other aspects that make László Lovász special. Mathematicians are
often divided into ``problem solvers'' and ``theory builders''.  Graph theory
is, in particular, a field to which problem solvers are drawn. Theory builders
often see deep and unusual connections, but often leave the difficult
exploration of details to others. As we will demonstrate, Lovász is a member of
this rare breed of people who possess both talents. Moreover, he brought his
talents to bear not only in one field of mathematics, he has also fertilized
and inspired significant developments in a wide range of other areas. If asked
to formulate the essence of his contributions in few words, we could use the
following three:

\begingroup
  \setlength{\pltopsep}{1\baselineskip}

\begin{compactdesc}
\item[Depth]Lovász solves many important and widely known problems in a
  competitive environment. He isolates seemingly special topics and develops
  them into broad and important calculi.

\item[Elegance] His solutions are often surprisingly (and sometimes
  seemingly) simple. At the same time, they often are mathematically beautiful
  and suggest fundamentally new ways to address a problem.

\item[Inspiration] Many of his solutions are the basis of further active
  research and even the foundations of whole new areas.
\end{compactdesc}
László Lovász published eleven books and more than 300~articles. There is no
way to survey his contributions in an article like this. We have chosen to
sketch some of the publications and topics that we consider highlights, are not
too difficult to explain, had significant impact, moved the frontier of
knowledge in the interface of mathematics and computer science substantially,
and are of lasting value.

\endgroup

\section{Logic and Universal Algebra -- Homomorphisms and Tarski's Problem}
\label{sec:GN2}

\begin{discussedreferences}
L. Lovász. Operations with structures. \emph{Acta Mathematica Academiae
Scientiarum Hungarica} 18:321--328, 1967.

L. Lovász. On the cancellation law among finite relational structures.
\emph{Periodica Mathematica Hungarica} 1:145--156, 1971.

M. Freedman, L. Lovász, L. Schrijver. Reflection positivity, rank
connectivity, and homomorphisms of graphs. \emph{Journal of American
Mathematical Society} 20(1):37--51, 2007.
\end{discussedreferences}

\noindent
Up to the 1960s graph theory was mainly concerned with graphs as objects. Graph
parameters were introduced and the structural properties of graphs having these
properties were investigated. László Lovász made, as we will outline, very
significant contributions to this kind of research, but he left his first
fundamental mark, when he was 19~years old, in the more general context of
universal algebra.

Intending to step out of the object orientation of graph theory, Lovász got
interested in operations with graphs and their algebraic properties.  We all
know that, for nonzero real numbers $a$, $b$, and $c$, the equation $ac = bc$
implies $a=b$. Suppose we have three graphs $A$, $B$, and $C$, and suppose we
have defined a product ``$\times$'' for which $A\times C=B\times C$ holds, can
we infer that $A=B$? Such a question makes only sense if equality ``$=$'' is
replaced by ``isomorphic'' and the concrete issue to be addressed is: Under
what conditions, does such a ``cancellation law'' hold?

Questions of this type were asked by Alfred Tarski, in the context of finite
relational structures, to students in Berkeley in the 1960s.  Lovász points
this out in the following quote, extracted from his article~\cite{Lovasz1967},
where he states the question and announces his solution:

\begin{figure}[H]
\centering
\includegraphics[width=\linewidth]{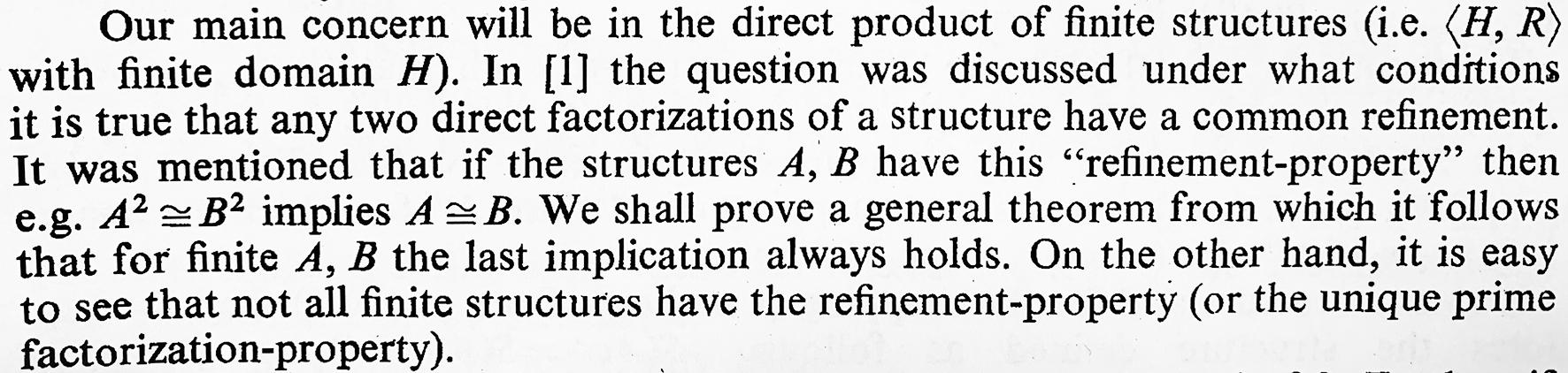}

\caption{Quote from~\cite{Lovasz1967}}
\label{fig:GN3}
\end{figure}

Reference [1] in the quote above is the paper~\cite{ChangJonssonTarski1964} of
Chang, Jónsson, and Tarski of 1964, see also~\cite{McKenzie1971}.

A finite \emph{graph} $G$ with \emph{vertex set} $V(G)$ and \emph{edge set}
$E(G)$ is such a relational structure where $V(G)$ is the ground set and the
edges $uv$ define the (binary) relations between vertices $u$ and $v$. A
standard product in graph theory is the \emph{direct} (also named categorial or
tensor) \emph{product} $G_{1}\times G_{2}$ of two graphs $G_{1}$ and $G_{2}$.
Its vertex set is $V(G_{1})\times V(G_{2}) = \{(u,v) \mid u\in V(G_{1}), v\in
V(G_{2})\}$ and its edge set $E(G_{1}\times G_{2})$ is defined to be the set of
all pairs of vertices $(u_{1},u_{2}),(v_{1},v_{2}) \in V(G_{1})\times V(G_{2})$
with $u_{1}v_{1} \in E(G_{1})$ and $u_{2}v_{2} \in E(G_{2})$. The question to
be addressed is: Given two graphs $G$ and $H$ and a third graph $F$, can one
conclude that $G$ and $H$ are isomorphic if the direct product $F\times G$ is
isomorphic to $F\times H$? This particular question and most of the related
problems for finite relational structures were unsolved, despite considerable
effort. The earlier solution approaches taken were usually elementary, trying
to reduce the problem to known invariants.

Lovász devoted to these problems three of his early papers written in~1967,
1971, 1972. His approach was radically different: He invented a new invariant
which solved these problems for the direct product in full generality. His
results completely changed this area.

The Lovász argument is easy and can be given here in full.  Interestingly,
young Lovász formulates his results very generally for finite \emph{relational
structures}, i.e., objects of the form $\pmb{A}  =  (X_{\pmb{A}},
(R_{\pmb{A}};R  \in  L))$ where $R_{A}$ is a subset of ${X\ }^{a(R)}$ ($a(R)$ is
the arity of the relational symbol $R$; $L$ is the fixed set of symbols usually
called language). Shortly, we speak about $L$-\emph{structures}.

A \emph{homomorphism} $f: \pmb{A} \rightarrow  \pmb{B} =
(X_{\pmb{B}}, (R_{\pmb{B}}; \pmb{R} \in  L))$
is a mapping 
$f: X_{\pmb{A}} \rightarrow X_{\pmb{B}}$ 
such that for every $R \in L$ holds
$(x_{1}, \ldots, x_{a(R)}) \in {R}_{\pmb{A}} \Rightarrow (f(x_1), \ldots, f  x_{a(R)}  \in R_{\pmb{B}}$.
The \emph{product} $\pmb{A} \times \pmb{B} $ is defined
as 
$X_{{\pmb{A}} \times \pmb{B}} = X_{\pmb{A}}  \times  X_{\pmb{B}}$
where 
$\pmb{R}_{{\pmb{A}} \times \pmb{B}}$ 
is the set of all tuples
$((x_{1},y_{1}), \ldots, (x_{a(R)},y_{a(R)}))$ 
where
$(x_{1},\ldots x_{a(R)}) \in R_{\pmb{A}}$ 
and
$(y_{1},\ldots,y_{a(R)}) \in R_{\pmb{B}}$.

Note that the projections 
$\pi_{\pmb{A}}: X_{\pmb{A} \times \pmb{B}} \rightarrow X_{\pmb{A}}$ 
and
$\pi_{\pmb{B}}: X_{\pmb{A} \times \pmb{B}} \rightarrow X_{\pmb{A}}$ 
are homomorphisms. Up to an isomorphism, projections uniquely determine the
above product.  (The whole theory may be restated in categorical terms as
worked out in papers by Lovász~\cite{Lovasz1972c} and Pultr~\cite{Pultr1973}.)

Denote by $\hom(\pmb{A},\pmb{B})$ the number of homomorphisms from
$\pmb{A}$ to $\pmb{B}$. The key of Lovász's argument is the following
statement:

\begin{nntheorem}
Finite $L$-structures $\pmb{A}$ and $\pmb{B}$ are isomorphic if and
only if for every other finite structure $\pmb{C}$ holds: 
$\hom(\pmb{C},\pmb{A}) = \hom(\pmb{C}, \pmb{B})$.
\end{nntheorem}

In other words (and today's setting), if we take a fixed enumeration
$F_1,F_2,\ldots F_n, \ldots$ of all non-isomorphic finite
graphs then the vector
$L(\pmb{A}) = (\hom(F_i, \pmb{A}); \ i = 1, \ldots)$ is
the isomorphism invariant, expressed equivalently:
$\pmb{A} \cong \pmb{B}$ if and only if
$L(\pmb{A}) = L(\pmb{B})$.

Hell and Nešetřil~\cite{HellNesetril2004} (and others) call this
invariant $L(\pmb{A})$ \emph{Lovász vector}.

This setting is very suitable for the Tarski problem. For example, one
immediately obtains that for finite structures $\pmb{A}^{k} \cong
\pmb{B}^{\mathbf{k}}$ holds if and only if $\pmb{A}\  \cong  \pmb{B}$. This
follows readily from $\hom(\pmb{C},\pmb{A}^k) = (\hom(\pmb{C},\pmb{A}))^k$.

For brevity we mention another consequence for the special case of graphs.
If~$\pmb{C}$ is a nonbipartite graph then $\pmb{A} \times \pmb{C} \cong
\pmb{B} \times  \pmb{C}$ if and only if $\pmb{A} \cong \pmb{B}$. (Note that for
bipartite graphs $\pmb{C}$, cancelation need not hold as already for
circuits we have ${2\pmb{C}}_3 \times K_2 \cong \pmb{C}_6 \times K_2$.)

The above theorem is very general and yet the proof is easy. In the
nontrivial direction we prove by induction on the cardinality
$|X_{\pmb{C}}|$ of the ground set
$X_{\pmb{C}}$ that, if $\hom(\pmb{C},\pmb{A}) =
\hom(\pmb{C},\pmb{B})$, then also the number of
injective homomorphisms coincides, i.e., $\inj(\pmb{C},\pmb{A}) = \inj(\pmb{C},\pmb{B})$.

In the inductive step we have $\hom(\pmb{C}, \pmb{A}) = \sum_{\theta}
\inj(\pmb{C}/\theta, \pmb{A})$ where $\theta$ is an equivalence on
$X_{\pmb{C}}$.  Thus by induction assumption we have $0 = \hom(\pmb{C},\pmb{A})
- \hom(\pmb{C},\pmb{B})  =  \inj(\pmb{A},\pmb{A}) - \inj(\pmb{A},\pmb{B}) =
\inj(\pmb{B},\pmb{B}) - \inj(\pmb{B},\pmb{A})$.  But obviously $\inj
(\pmb{A},\pmb{A}) > 0$ and $\inj(\pmb{B},\pmb{B}) > 0$, and thus, we have that
there are injective homomorphisms from $\pmb{A}$ to $\pmb{B}$ and also from
$\pmb{B}$ to $\pmb{A}$. Now as $\pmb{A}$ and $\pmb{B}$ are finite structures we
have that $\pmb{A} \cong \pmb{B}$.

Lovász recognized in the Tarski problem a magnificent pearl. His theorem turned
out to be very useful. It found many applications and inspired further
research. This continues until today, see the articles by Lovász and
Schrijver~\cite{LovaszSchrijver2009}, Dvořák~\cite{Dvorak2009} and Dawar et
al.~\cite{DawarJaklReggio2021}, for example.

The papers~\cite{Lovasz1967} and~\cite{Lovasz1971} of Lovász belong to the
first occurrences of homomorphisms in graph theory. Their successful
utilization led to a rich calculus (see, e.g., the books by Hell and
Nešetřil~\cite{HellNesetril2004} and Lovász~\cite{Lovasz2012} and the article
by Borgs et al.~\cite{BorgsChayesLovaszSosVesztergombi2006}). We outline
important parts of this approach.

Lovász already defined in~\cite{Lovasz1967} exponential structures
$\pmb{A}^{\pmb{B}}$ (and exponential 
\linebreak
graphs~$G^{H}$). These played very
recently a decisive role in the disproof of the Hedetniemi conjecture which
claimed that $\chi(G\times H)  = \min(\chi(G), \chi(H))$, see
Shitov~\cite{Shitov2019}, Wrochna~\cite{Wrochna2019}, Tardif~\cite{Tardif2022},
and Zhu~\cite{Zhu2021}.

Another application of homomorphism counting was provided by Lovász
in~\cite{Lovasz1972d} which deals with the following problem: When can one
recognize a given finite structure from the collection of all its proper
substructures? The special case for undirected graphs is a classical
\emph{conjecture of Ulam}, see~\cite{Ulam1969}, which may be formulated in our
setting as follows:

\begin{noheadingthm}
Do the homomorphism numbers $\hom(F,G)$ for all graphs $F$ with fewer edges
than~$G$ determine the graph $G$?
\end{noheadingthm}

This conjecture is known to be true for special classes of graphs (such
as trees and maximal planar graphs), and the proofs usually consist of a
complicated case analysis. In~\cite{Lovasz1972d} Lovász gave the first
general result: The conjecture is true for graphs that have more edges
than their complement (i.e., more than half of all edges).

The proof, although not directly linked to the above theorem proceeds
again by clever homomorphism counting. Shortly after, this proof has
been extended by Müller in~\cite{Mueller1977} (again by homomorphism
counting) to graphs with $n \log n$ edges. This is still the
best result.

Counting of homomorphisms and the investigation of their structure are
cornerstones of further areas of mathematics and theoretical computer science.
We just indicate three examples, where they play important roles: Tutte
polynomials and their variants, see~\cite{EllisMonaghanMoffatt2022};
constraint satisfaction problems (which can alternatively be viewed as
existence theorems for general relational structures),
see~\cite{FederVardi1998}; and partition functions in statistical physics,
see~\cite{CaiChen2017, Lovasz2012, BorgsChayesLovaszSosVesztergombi2012}.

Let us finally elaborate on partition functions, the last item mentioned
above. The concept of graph homomorphisms can be extended to graphs with
loops and weights assigned to vertices $\propto_{v} (G)$ and
edges ${\beta}_{uv}(G)$. For unlabelled graphs $F$ and labelled
graphs $G$, one can define naturally the weight of a mapping
$\varphi: V(F) \rightarrow V(G)$ and then the total weight of
$\hom(F,G)$. Allowing weights on the vertices and edges greatly extends
the expressive power of (weighted) homomorphisms. For example, the
number $\hom(F,G)$ can express the number of colorings (leading to
chromatic and Tutte polynomials), the counting of stable sets
(corresponding to the so-called \emph{hard core model} in statistical
physics) and also the counting of nowhere zero flows and $B$-flows
(i.e., flows attaining values from a given set $B$ only). All these
are parameters of the form $\hom(-,H)$. Freedman, Lovász, and
Schrijver~\cite{FreedmanLovaszSchrijver2007} provided a structural
characterization for all such parameters as follows:

\begin{nntheorem}
Let $f$ be a (real) graph parameter defined on
multigraphs without loops. Then f is equal to $\hom(-,H)$
for some weighted graph H on q vertices if and only if
$f(K_{0}) = 1$, the $f$ connection matrix $M(f,k)$ is
reflection positive, and its rank satisfies $r(M(f,k)) \leq q^{k}$
for all $k \geq 0$.
\end{nntheorem}

(Briefly: Above, $K_k$ is the complete graph on $k$ vertices;
the connection matrix $M(f,k)$ is defined by values of the parameter
$f$ for amalgams of $k$-multilabeled multigraphs; reflection
positivity means that, for all $k$, such matrices are positive
semidefinite.)

This theorem led to many similar results for other classes of graphs and for
other types of homomorphism numbers (e.g., in a dual setting with $\hom(F,-)$
instead of $\hom(-,H)$, see~\cite{LovaszSchrijver2010}). In terms of
statistical physics, this theorem can be viewed as a characterization of
partition functions of vertex coloring models.

Lovász wrote extensively on this topic and devoted -- ten years ago -- a
monograph~\cite{Lovasz2012} to this subject, where the topics indicated
here are treated in depth.

\section{Coloring Graphs Constructively (on a Way to Expanders)}
\label{sec:GN3}

\begin{discussedreferences}
L. Lovász. On chromatic number of finite set-systems. \emph{Acta
Mathematica Academiae Scientiarum Hungaricae}, 19:59--67,1968.
\end{discussedreferences}

\noindent
The \emph{chromatic number} $\chi(G)$ of a graph $G$ is the minimum
number of colors which suffice to color all vertices of $G$ such that no
two adjacent vertices get the same colour. Alternatively, using the
notion of the preceding section, $\chi(G)$ is smallest $k$ for
which $\hom(G, K_{k}) > 0$.

The chromatic number belongs to the most frequently studied
combinatorial parameters. Reasons for such an attention are that the
question of how to color the countries on a map can be easily explained
to everyone and that the mathematical modelling of this question can be
employed as an appealing introduction to graph theory. The ``colorful
story of the 4-color conjecture'' can be used to shed some light on the
rich history of mathematics and the difficulty of finding proofs for
problems that appear to be easy. Coloring the vertices of a graph
captures the substance and the difficulty of many problems. In a
multiple sense, the chromatic number is a difficult concept.

Just consider the easiest question: Are there graphs with large
chromatic number? Of course, complete graphs $K_{n}$ satisfy $\chi$
($K_{n})= n$. But are there any other essentially different
graphs?

The answer is yes and a classical result, rediscovered several times,
states that, for every $k\geq1$, there are graphs $G_k$ for
which $\chi(G_{k}) = k$ and $G_{k}$ does not contain
$K_3$ (i.e., the triangle) as a subgraph. This result and its many
ramifications, for instance in extremal graph theory, are still in the
current focus of coloring research. In fact, any new constructive proof
of the existence of such graphs $G_{k}$ is interesting and attracts
great attention. Here is perhaps the simplest proof of this fact: Let us
define, for any integer $n\geq4$, the graph $G = (V,E)$ where $V$
is the set of integer pairs
$\{ ij\}, 1 \leq i  < j \leq n$, and
$\{ij, kl\}  \in E \text{ if }  i < j = k < l$. Such a graph
$G\ $is called a \emph{shift graph}. $G$ has no triangles, and it can be
shown that $\chi(G) = \lbrack \log n \rbrack$.

But this is not the end of the story. Graphs may have high chromatic
number and very low edge density. P.~Erdős showed in~\cite{Erdos1959}
that there exist graphs which have arbitrarily large chromatic number
and which are locally trees and forests.

\begin{nntheorem}
For every $k, l$ there exists a graph ${G}_{k,l}$ such that
$\chi({G}_{k,l})\geq k$ and ${G}_{k,l}$ does not contain circuits of length
$\leq l$. (So the shift graph above is a graph of type~$G_{k,3}$.)
\end{nntheorem}

Erdős' proof was a landmark. It constitutes one of the key applications of the
probabilistic method in graph theory, see, e.g.,~\cite{AlonSpencer2016}. The
proof shows that the probability of the existence of such graphs $G_{k,l}$ is
positive, but does not give any hint how to construct concrete examples of
graphs of type~$G_{k,l}$. The construction of such graphs has been a
longstanding problem with very slow progress (for the historic development and
related issues, see, e.g., the Nešetřil article~\cite{Nesetril2013}).

The first constructive proof of the theorem above was found by Lovász in one of
his early papers~\cite{Lovasz1968}. It was one of the highlights of the 1969
conference in Calgary; and through his proof, Lovász again changed the setting
of the problem as he constructed the graphs $G_{k,l}$ as special cases of a
more general theorem about hypergraphs. His complicated construction was later
simplified, the Nešetřil--Rödl construction is perhaps the
simplest~\cite{NesetrilRoedl1979}.

But various problems remained.
  
One of them is the question whether one can provide a construction that
uses only graphs. The answer is positive. I.~Kříž~\cite{Kriz1989} and
more recently N.~Alon et al.~\cite{AlonKostochkaReinigerWestZhu2016} came
up with such constructions, and Ramanujan graphs have to be mentioned
here as well.

The existence of graphs $G_{k,l}$ with a large chromatic number and
no short circuit is a phenomenon of finite (and of countable) graphs.
For graphs with an uncountable number of vertices and uncountable
chromatic number, an analogous result does not hold. This was shown by
Erdős and Hajnal~\cite{ErdosHajnal1966}:

\begin{noheadingthm}
If the chromatic number of a graph is uncountable then it contains
every bipartite graph.
\end{noheadingthm}

A consequence of this result is that such a graph contains every circuit
of even length, for example the circuit $\pmb{C}_{4}$ of length four.

Graphs $G_{k,l}$ are what can be called difficult examples. They also
play an important role in Ramsey theory, extremal combinatorics,
topological dynamics, and model theory, to name just a few. In all these
areas they are used as examples of complex yet locally simple
structures; they are prototypes of local-global phenomena.

It took some time to understand why the construction of graphs $G_{k,l}$
matters, why it is important to know such graphs explicitly. This led to an
explosion of theoretical developments combining group theory, number theory,
geometry, algebraic graph theory, and, of course, combinatorics. The key
notions are now familiar to every student of theoretical computer science:
expanders, Ramanujan graphs and sparsification, see
Margulis~\cite{Margulis1973}, Lubotzky, Phillips, and
Sarnak~\cite{LubotzkyPhillipsSarnak1988}, and Spielman and
Teng~\cite{SpielmanTeng2011}.

An expander graph, for instance, is a finite, undirected multigraph
(parallel edges are allowed) in which every subset of the vertices that
is not ``too large'' has a ``large boundary''. There are various
formalizations of these notions. Each of them gives rise to a different
notion of expanders, e.g., edge expanders, vertex expanders, and
spectral expanders. Expander graphs have found applications in the
design of algorithms, error correcting codes, pseudorandom generators,
sorting networks, robust computer networks and hash functions in
cryptography. They also played a role in proofs of important results in
computational complexity theory, such as the PCP theorem.

The construction and structure of graphs similar to $G_{k,l}$ continues to be
one of the key problems of finite combinatorics and has a character of a saga
(see, e.g., Hoory et al.~\cite{HooryLinialWigderson2006} and
Nešetřil~\cite{Nesetril2013}).

Coloring of graphs and hypergraphs has been a permanent theme of Lovász,
and thus, it is mentioned in most sections of our survey. For example,
one of the motivations of the next section was the study of 3-chromatic
linear hypergraphs, i.e., hypergraphs in which edges meet in at most
one vertex, or equivalently, hypergraphs without cycles of length~2.

\section{The Lovász Local Lemma}
\label{sec:GN4}

\begin{discussedreferences}
P. Erdős, L. Lovász. Problems and results on 3-chromatic hypergraphs and
some related questions. In \emph{Infinite and Finite Sets. Coll. Math.
Soc. J. Bolyai}, North Holland:609--627, 1975.
\end{discussedreferences}

\noindent
A \emph{hypergraph} is a collection of sets. The sets are called
\emph{edges}, the elements of the edges are \emph{vertices}. The
\emph{degree} of a vertex is the number of edges containing it. A
hypergraph is called \emph{r-uniform} if every edge has $r$~vertices. The
\emph{chromatic number of a hypergraph} is the least number $k$ such that
the vertices can be $k$-colored so that no edge is monochromatic.

\begin{figure}
\centering
\includegraphics[width=.8\linewidth]{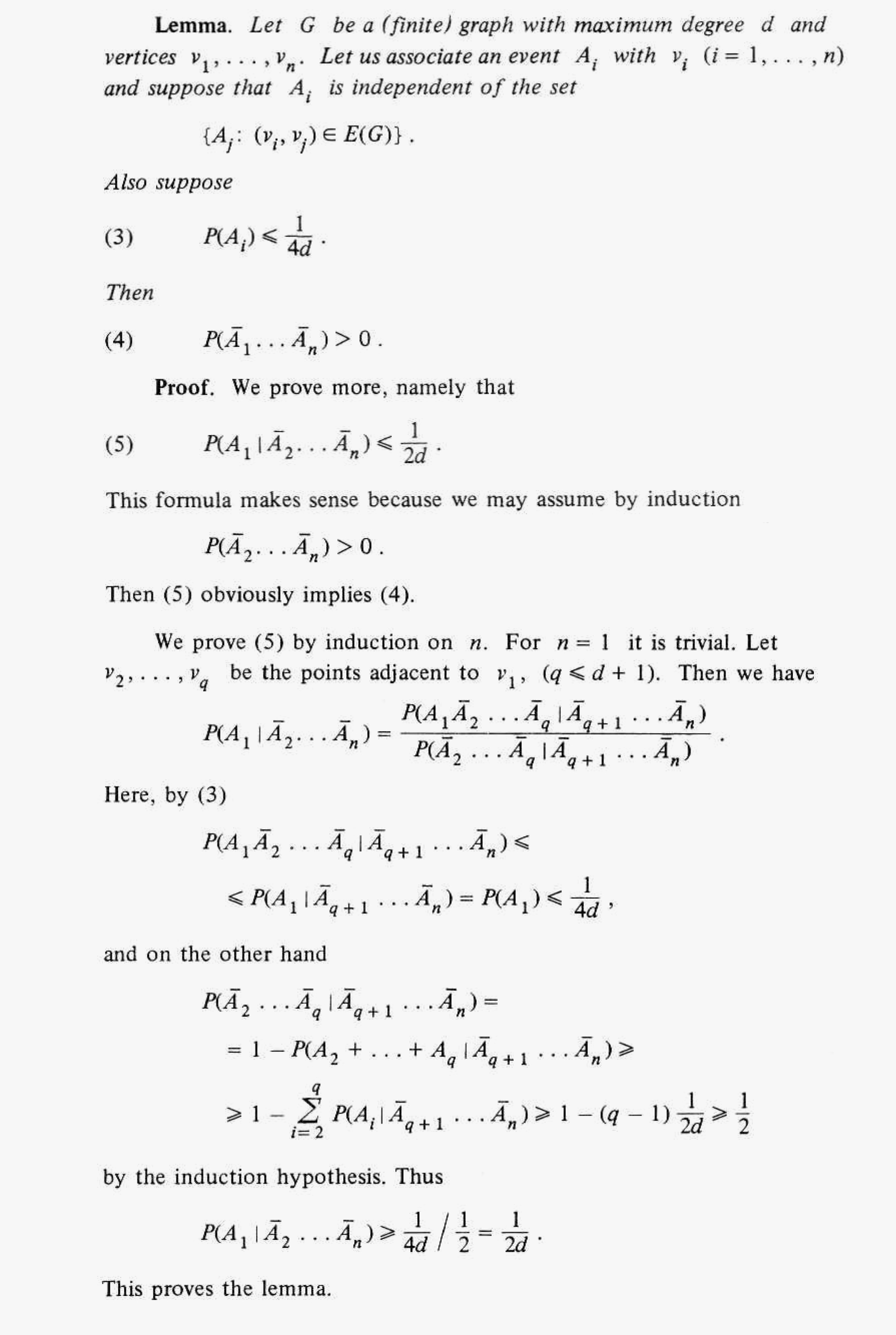}

\caption{Extracted from~\cite{ErdosLovasz1975}}
\label{fig:GN4}
\end{figure}

Graphs with chromatic number at least 3 are simple to characterize: they must
contain an odd circuit. But for hypergraphs, even the characterization of
3-chromatic 3-uniform hypergraphs is difficult (it is an \cNP-complete
problem). Lovász and Woodall had independently shown that every 3-chromatic
$r$-uniform hypergraph contains a vertex of degree at least $r$. Erdős and
Lovász~\cite{ErdosLovasz1975} aimed at generalizing this result in various
ways. One of the key results of their article is the following:

\begin{nntheorem}
A $(k+1)$-chromatic r-uniform hypergraph contains an edge which is
intersected by at least $k^{r-1}/4$ other edges. Thus, the
degree of at least one vertex is larger than $k^{r-1}/(4r)$.
\end{nntheorem}

To prove this theorem, the authors employed probability theory. As
pointed out by Erdős, Lovász contributed to the proof a substantial new
result of elementary probability. This was later called the \emph{Lovász
Local Lemma}.

The motivation for this lemma comes from a well-known observation of
elementary probability:

If $X_{1}, \ldots,X_{n}$ are random events which are pairwise
independent and if the probability of each event $X_{i}$ is
smaller than $1$, then the probability that none of the events
$X_{i}$ occurs is positive. The Lovász Local Lemma is a quantitative
refinement of this observation for variables which are dependent.

Fig.~\ref{fig:GN4} shows the formulation of the Lovász Local Lemma as stated
and proved in the original article~\cite{ErdosLovasz1975}. Indeed, it is ``just
a lemma''.

Crystal clear: Not only when the events are independent, but if the
dependence graph $G$ has a small degree ($\leq d$) then also none of
the events occurs with positive probability. The adjective local in the
name of the lemma refers to the situation that each event is dependent
only on a small number $d$ of others.

It is hard to overestimate the general importance of this result that just
turned up as a ``supporting observation'' for a proof in the chromatic theory
of hypergraphs. It appears again and again in multiple applications,
ramifications, and forms. It is not possible to cover here all the applications
in Ramsey theory (see Spencer~\cite{Spencer1975}), extremal combinatorics (see
Alon and Spencer~\cite{AlonSpencer2016}), number theory, and elsewhere (see,
e.g., Ambainis et al.~\cite{AmbainisKempeSattath2010}, He et
al.~\cite{HeLiSunZhang2019}, and Szegedy~\cite{Szegedy2013}). It was also
discovered, see~\cite{ScottSokal2005}, that the Lovász Local Lemma closely
relates to important results of Dobrushin in statistical
physics~\cite{Dobrushin1996}. In fact, the proper setting of the Dobrushin
results is in the context of graph limits, see~\cite{Lovasz2012}, which we
discuss in Section~\ref{sec:GN17}.

One of the motivations for~\cite{ErdosLovasz1975} is the following number
theoretic problem which goes back to Ernst Straus (who was an assistant
of Albert Einstein): Is there a function~$f(k)$ such that, if $S$ is any set
of integers with $|S| = f(k)$, then the integers can be
$k$-colored so that each color meets every translated copy of $S$ (i.e.,
every set of the form $S + a = \{x+a \mid a \in S\}$)? Lovász and
Erdős, already in their paper~\cite{ErdosLovasz1975}, made use of the
Lovász Local Lemma to prove the following more geometric generalization
of the question asked by Straus:

\begin{noheadingthm}
For every $k$, there exists a function $f(k)$, such that $f(k) \leq k \log k$
and for every set~$S$ of lattice points in the $n$-dimensional space $E^{n}$
with $|S| > f(k)$ there exists a $k$-coloring of all lattice
points such that each translated copy of $S$ contains points of all $k$ colors.
\end{noheadingthm}

A side remark: There are many variants of coloring problems, and some of
them are surprisingly difficult. For example, during a conference in
Boulder in~1972 Paul Erdős, Vance Faber, and László Lovász asked whether
the vertices of any $n$-uniform linear hypergraph with $n$ edges can be colored
by $n$ colors such that the vertices of any edge get all $n$ colors. This
question has many reformulations and turned out to be more difficult
than originally thought (even by the authors as Erdős originally offered
\$50 for a solution and eventually increased the prize to \$500). About
50~years later the Erdős--Faber--Lovász conjecture was shown to be true
for large values of $n$ by D.\,Y.~Kang, T.~Kelly, D.~Kühn, A.~Methuku, and
D.~Osthus~\cite{KangKellyKuehnMethukuOsthus2023}.

Nowadays, the Lovász Local Lemma is a ``standard trick'' which is often
taught in basic courses. And it is a very effective trick, as Joel
Spencer once remarked: ``\emph{Using the Local Lemma one can prove the
existence of a needle in a haystack}.''

But the Lovász Local Lemma delivers only existence. The above proof does not
yield a method how to find that needle. We only know that certain things exist
with positive probability. Only much later a constructive proof was found by
Moser and Tardos~\cite{MoserTardos2010}. (Remark: A constructive proof for
the above Straus' problem is in Alon et al.~\cite{AlonKrizNesetril1995}; see
also J.~Beck~\cite{Beck1991}.) Recently, Harvey and
Vondrák~\cite{HarveyVondrak2020} found another constructive approach to the
Lovász Local Lemma.

Investigations of infinite versions (Borel and measurable) of the Lovász Local
Lemma started also very recently by A.~Bershteyn, G.~Kun, O.~Pikhurko, and
others, see, e.g.,~\cite{Bernshteyn2019}.)

The Lovász Local Lemma became what one can truly call \emph{a combinatorial
principle}. This is László Lovász at its best: Maybe no other
Lovász-contribution is so profoundly simple and yet useful and elegant.

\section{Coloring Graphs via Topology}
\label{sec:GN5}

\begin{discussedreferences}
L. Lovász. Kneser's conjecture, chromatic number, and homotopy.
\emph{Journal of Combinatorial Theory A} 25:319--324, 1978.
\end{discussedreferences}

\noindent
Combinatorial questions are often easy to formulate; some have also an
elementary solution. But in many cases, the elementary nature of
combinatorial problems is just the top of an iceberg, and the hidden
complexity must be discovered and tamed before a solution can be found.

A beautiful example of this is the following elementary problem posed in
1955 by Martin Kneser~\cite{Kneser1955} who was working on quadratic
forms. In today's language:

\begin{noheadingthm}
Let $X$ be a set with $n$ elements, $n \geq
2k > 0$. Denote by $\left(\begin{smallmatrix}X\\k\end{smallmatrix}\right)$
the set of all $k$-element subsets
of $X$. Then, for every coloring of the sets in
$\left(\begin{smallmatrix} X\\k\end{smallmatrix}\right)$ by fewer than 
$n - 2k + 2$ colors, there are two disjoint sets of the same color.
\end{noheadingthm}

\begin{figure}
\centering

\newcount\nodecount
\tikzgraphsset{
  declare={subgraph N}%
  {
    [/utils/exec={\global\nodecount=0}]
    \foreach \nodetext in \tikzgraphV
    {  [/utils/exec={\global\advance\nodecount by1}, 
      parse/.expand once={\the\nodecount/\nodetext}] }
  },
  declare={subgraph C}%
  {
    [cycle, /utils/exec={\global\nodecount=0}]
    \foreach \nodetext in \tikzgraphV
    {  [/utils/exec={\global\advance\nodecount by1}, 
      parse/.expand once={\the\nodecount/\nodetext}] }
  }
}

\begin{tikzpicture}[every node/.style={draw,circle,thin}, font=\footnotesize]
  \graph [clockwise,math nodes] {     
    subgraph C [V={ {$\{1,2\}$}, {$\{4,5\}$}, {$\{1,3\}$}, {$\{2,5\}$}, {$\{3,4\}$} }, name=A, radius=2.5cm]; 
    subgraph N [V={ {$\{3,5\}$}, {$\{2,3\}$}, {$\{2,4\}$}, {$\{1,4\}$}, {$\{1,5\}$} }, name=B, radius=1.1cm];
    \foreach \i [evaluate={\j=int(mod(\i+1,5)+1);}] in {1,...,5}{
      A \i -- B \i;  
      B \i -- B \j;
    }
  }; 
\end{tikzpicture}

\caption{$KG(5.2)$ = The Petersen graph}
\label{fig:GN5}
\end{figure}
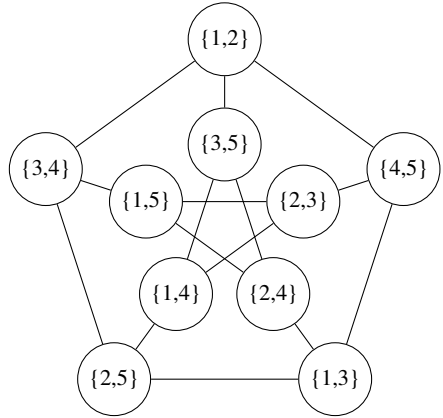

This problem can be reformulated as a graph theory question as follows.  Let
$KG(n,k)$ denote the graph (called \emph{Kneser graph}) whose vertices are all
$k$-element subsets of set $X=\{1,2,\ldots, n\}$, and in which two vertices are
joined by an edge if the corresponding $k$-element subsets are disjoint. For
example, $KG(n,1)$ is the complete graph $K_{n}$ and $KG(5,2)$ is the famous
Petersen graph (the ``universal'' counterexample to many conjectures in graph
theory) shown in Fig.~\ref{fig:GN5}.

Kneser's question reads now: \emph{Does the Kneser graph $KG(n,k)$
have chromatic number $n-2k+2$?}

It is easy to see that $\chi(KG(n,k)) \leq n - 2k + 2$.  However, to find the
fitting lower bound for the chromatic number proved to be much harder.

Lovász~\cite{Lovasz1968} solved this problem in a surprising way using
methods of algebraic topology. The general idea is the following. Lovász
associates with any graph $G$ a topological space and establishes a
connection between a topological invariant of this space with the
chromatic number of $G$. He then infers properties of the chromatic number
of $G$ from properties of the topological invariant of the associated
topological space. That this is possible and that topology can yield
solutions of difficult graph theory questions was completely unexpected.
Lovász's success with this approach was the starting point of a new
field: \emph{topological combinatorics}. We briefly sketch the main
steps of Lovász's solution of Kneser's problem here.

Lovász proceeds as follows: Given a graph $G = (V,E)$, the
neighborhood of a vertex~$v$ is composed of all vertices adjacent
to $v$ in $G$. The \emph{neighborhood complex}~$N(G)$ of
$G$ consist of all the vertices $V$ of the graph $G$;
the simplices of $N(G)$ are sets of vertices with a common
neighbor in the graph. Homomorphisms between graphs lead to
continuous mappings of neighborhood complexes. From the topological
connectivity of $N(KG(n,k))$ it is possible to construct an
antipodal continuous mapping between spheres
$(N(K_{m+2})$ is an \emph{m}-dimensional sphere) and
one can then apply the Borsuk--Ulam theorem. Thus, Lovász obtained:

\begin{nntheorem}
If the neighborhood complex $N(G)$ of a graph G is (topologically)
$k$-connected then $\chi(G) \geq k + 3$.
\end{nntheorem}

(Topologically $k$-connected means that there are no holes of
dimension $\leq k$. For (simply) connected complexes this is
equivalent to the fact that all $i$-homological groups vanish for
$i = 0,1, \ldots, k$.)

Lovász finally proves a theorem on the connectivity of neighborhood
complexes of graphs from which he can infer that the neighborhood
complex of a Kneser graph $N(KG(n,k))$ is topologically
$(n-2k-1)$-connected. This establishes that the Kneser graph $KG(n,k)$ has
chromatic number $n-2k+2$.

This connection (and the whole proof) immediately led to intensive research.
Other proofs of this theorem were found (among them ``book proofs'' of
Barany~\cite{Barany1978} and Green~\cite{Green2002}), but all lower bounds for
the chromatic number of Kneser graphs use or at least imitate Lovász's
topological proof. Matoušek's book~\cite{Matousek2003} surveys in detail
various implications and modifications of the proof techniques. For example, it
has been shown in~\cite{GyarfasJensenStiebitz2004} that the $k$-times
generalized Mycielski construction has chromatic number $k+2$, and again,
topological arguments are the basis of the only known proof of this fact. The
paper~\cite{GyarfasJensenStiebitz2004} contains the following interesting
construction of graphs~$G_{k}$.

Put $[k] = \{1,2,\ldots,k\}$. The vertices of
$G_{k}$ are all pairs $(i,A)$ where $i
\notin A$ and $A$ is a nonepty subset of $\lbrack k\rbrack$.
$(i,A)$ and $(j,B)$ form an edge in $G_k$ if $i \in B, j \in A$ and $A$ and $B$ are
disjoint. This ``Kneser-like'' graph $G_{k}$ has
remarkable properties: Its chromatic number is $k$, it is critical
(i.e., every proper subgraph has a smaller chromatic number) and every
strongly $k$-colorable graph has a homomorphism into it; it is the
unique graph with this property. (A \emph{strong coloring} of a graphs
is a coloring where the
neighborhood of any color class forms a stable set. Such a graph
obviously has no triangles.)

The only known proof of these properties is an adaptation of Lovász's
topological proof.

These examples of graphs were instrumental in the recent disproof of the
Hedetniemi conjecture (that intended to establish a connection between
the direct product of two graphs and their chromatic number which we
mentioned in Section~\ref{sec:GN2}; see also~\cite{SimonyiZsban2010}) and also in the
study of \emph{gap problems} for constraint satisfaction problems.
Related questions in this area are called \emph{promised problems}. A
typical question here is: How difficult is to 5-color graphs or
hypergraphs under the assumption that we know they are 3-colourable, see~\cite{DinurRegevSmyth2005, BartoKozik2022}, and~\cite{Wrochna2022}.

Lovász's paper opened a whole area whose fruits are still continuing to
appear. Matoušek in the preface to~\cite{Matousek2003} rightly wrote that
Lovász's proof of the Kneser conjecture is a masterpiece of imagination.

\begin{figure}[t]
\centering
\includegraphics[width=.75\linewidth]{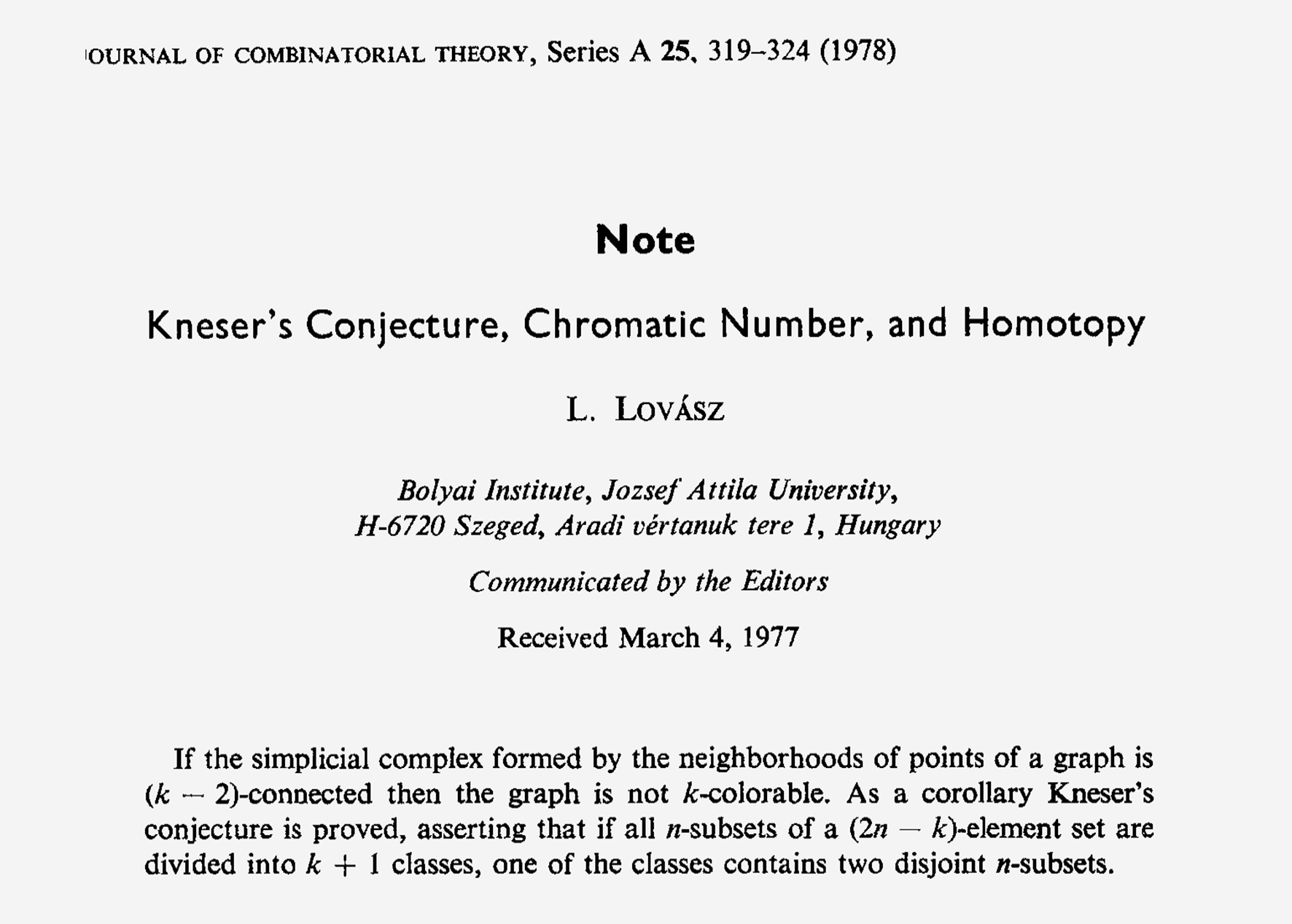}

\caption{The beginning of the article~\cite{Lovasz1978} starting topological
combinatorics}
\label{fig:GN6}
\end{figure}

Yet, in typical Lovász style, it was published just as a note (see Fig.~\ref{fig:GN6}).

Lovász's solution of the Kneser problem did not exhaust his topological
imagination nor the potential of topological methods in combinatorics. He
returned to this approach frequently during his career, often in collaboration
with Lex Schrijver. We mention just one of the highlights of their cooperation. 

Motivated by estimating the maximum multiplicity of the second eigenvalue of
Schrödinger operators, Colin de Verdi\`ere introduced a new invariant for graphs
$G$, denoted $\mu(G)$, based on spectral properties of matrices associated with
$G$. He proved that $\mu(G) \leq 1$ if and only if $G$ is a disjoint union of
paths, that $\mu(G) \leq 2$ if and only if $G$ is outerplanar, and that $\mu(G)
\leq 3$ if and only if $G$ is planar. 

Robertson, Seymour, and Thomas showed that a graph $G$ is linklessly embeddable
if and only if $G$ does not have any of the seven graphs in the Petersen family
as a minor. Their combinatorial result implies that $\mu(G) \leq 4$ if $G$ is
linklessly embeddable, and they conjectured that $\mu(G) \leq 4$ if and only if
G is linklessly embeddable. Lovász and Schrijver, see
\cite{LovaszSchrijver1998}, proved the only if part of this topological
characterization. The key ingredient of their proof is a new Borsuk-type
theorem on the existence of antipodal links, which is an extension of a
polyhedral version of Borsuk's theorem due to Bajmóczy and B\'ar\'any. The
combination of all these results provides a fascinating characterization of
graphs $G$ satisfying $\mu(G) \leq 4$ by means of spectral, combinatorial, and
topological properties. Topological methods seem to keep on flourishing in
combinatorics and graph theory.

\section{Geometric Graphs and Exterior Algebra}
\label{sec:GN6}

\begin{discussedreferences}
L. Lovász. Flats in matroids and geometric graphs. In
\emph{Combinatorial Surveys}. Proc. 6 British Comb. Conf. Academic
Press, pages 45--86, 1977.
\end{discussedreferences}

\noindent
Many of Lovász's proofs deal with graphs (and hypergraphs) and make use of some
additional structures. The Shannon Capacity paper, see Section~\ref{sec:GN8},
involved a geometric structure which was added (orthogonal representation) so
that the problem could be solved. To solve the Kneser problem, discussed in
Section~\ref{sec:GN5}, Lovász employed results from topology.  To recognize
that methodology from other mathematical fields can be utilized, needs of
course mathematical maturity, skill, and imagination.  We want to highlight
that this is a different strategy than merely studying embeddings of graphs
(e.g., graphs on surfaces): the special embeddings are being incorporated in
proofs as tools in order to solve a (different) problem.

A very special example of this is the Lovász-article \cite{Lovasz1977} which is
a remarkable paper for multiple reasons.

The paper was published as an invited lecture in the proceedings of
6th British Combinatorial conference. These
proceedings volumes usually contain surveys of recent developments. In
contrast, the Lovász paper -- full of new ideas -- solved an important
problem and unleashed research in two different areas: First, it started
research in graphs where vertices are forming a matroid; Lovász uses
here the term \emph{geometric (or pregeometric) graphs}, and this
generalization is essential for solving the problem. Secondly, the paper
started the application of exterior algebra in combinatorics.
Particularly, Lovász defined exterior calculus in matroids and Grassman
graded matroids.

Why is Lovász introducing this general machinery? Well, he is explicit
about that in the introduction: 

\begin{quote}
\emph{This paper was intended to deal with the covering problems in graphs. It
has turned out, however, that their study becomes much simpler if a more
general structure, which we shall call geometric graph, is considered}.
\end{quote}
Lovász later on used the term geometric graph in a broader sense, and he
recently wrote the book~\cite{Lovasz2019} treating the whole area in
detail.

What were the ``covering problems'' of~\cite{Lovasz1977}?

The starting point was an old problem due to Tibor Gallai related to
$\tau$-critical graphs: The \emph{covering number} of graph
$G = (V,E)$, usually denoted by $\tau(G)$, is the minimum
cardinality of a set $A \subseteq V$ such that every
edge of $G$ meets $A$. (Such a set~$A $ is also called
\emph{hitting set}.)

$\tau(G)$ is a ``hard'' combinatorial parameter (ultimately related to
the stability number $\alpha(G)$ and the chromatic number
$\chi(G)$).

One approach to gain information about the covering number is to
consider graphs that are critical with respect to this parameter. A
graph $G = (V,E)$ is $\tau$\emph{-critical} if
$\tau(G) > \tau(G - e)$ for every edge $e \in E$. Gallai proved
in 1961 that every $\tau$-critical graph $G$ satisfies
$|V| \leq 2\tau(G)$. So given $\tau$ there are only finitely many
$\tau$-critical graphs and this implies a ``finite basis theorem''.

However, a much stronger statement holds. Let us denote the gap in the
above inequality by $\delta(G): = 2\tau(G) - |V(G)|$.
Then one can observe that, given a $\tau$-critical graph $G,$ the
graph $G'$ obtained from $G\ $by subdividing an edge of $G$ by an
even number of vertices is also $\tau$-critical, and obviously
$\delta(G) = \delta(G')$. Gallai conjectured that this is the only
operation that does not destroy $\tau$-criticality and that the number
of $\tau$-critical graphs with a given value $\delta$ is
(essentially) finite. And this was the motivation of Lovász for his
paper~\cite{Lovasz1977} in which he proved this conjecture.

\begin{nntheorem}
The number of connected $\tau$-critical graphs $G$
with gap
$\delta(G) = 2\tau(G) - |V(G)|  = \delta$ and all vertex degrees
$\geq 3$ is at most $2^{{5\delta}^{2}}$.
\end{nntheorem}

The proof of this result is complex. In fact, Lovász develops several new
tools. The whole paper makes effective use of geometric graphs (where vertices
form a matroid). This allows Lovász to carry on a subtle refinement of
induction procedures. He makes magnificent use of his vast experience with
matchings and generalized factors (this was the subject of his doctoral thesis
supervised by T.~Gallai) which found its way into his early book on matching
theory~\cite{LovaszPlummer1986} with M.~Plummer. The proof also implicitly
contains the ``skew Bollobás theorem'' (in a matroid setting) about an extremal
problem for set intersections of pairs of sets and many other inspiring ideas,
in particular, the surprising utilization of exterior algebra. This aspect of
the paper~\cite{Lovasz1977} also generated a whole new theory.

We shall illustrate the use of exterior algebra by the simpler example
of the (Prague) dimension of graphs (treated in another Lovász paper~\cite{LovaszNesetrilPultr1980}).

It is easy to prove that every graph is an (induced) subgraph of the direct
product of complete graphs (the product we introduced in
Section~\ref{sec:GN2}). The smallest number of such a set of complete graphs is
called \emph{dimension} $\dim(G)$ of the graph $G$.

Thus, $\dim(K_{n}) = 1$ and $\dim( K_{n} \times K_{n} \times \cdots \times 
K_{n}) \leq t$ (direct product of $t$ copies of~$K_{n}$).

It is very nice that we have equality here. The proof
in~\cite{LovaszNesetrilPultr1980} is one of the first applications of exterior
algebra in combinatorics which was initiated in~\cite{Lovasz1977}.

\begin{nntheorem}
$\dim(K_{n}^{t}) = t$ for every $t \geq 1,\ n \geq 2$.
\end{nntheorem}

$K_{2}^{2}$ is isomorphic to $K_{2} + K_{2}$ and $K_{2}^{t}$ is
isomorphic to a perfect matching (i.e., disjoint edges) of size
$2^{t - 1}$.

It suffices to prove $\dim( K_{2}^{t}) \geq t$. Given a
representation $f:K_{2\ }^{t} \rightarrow K_{N}^{d},$ we put
explicitly:

$f(i) = a_{i} = (a_{i}^{1},\ldots,a_{i}^{d})$ and
$f(i') = b_{i} = (b_{i}^{1},\ldots,b_{i}^{d})$ (we
think of matchings having edges
$\{ i,\ i'\}\ i = 1,\ldots, 2^{t - 1}$). Clearly all
these $2^{t}$ vectors are distinct.

The condition that $f$ is an embedding can be then captured by \smash{$\prod_{k
= 1}^{d}{(a_{i}^{k} - b_{j}^{k}) \neq 0}$} if and only if $i = j, \prod_{k =
1}^{d}{(a_{i}^{k} - a_{j}^{k}) = 0}$, and $\prod_{k = 1}^{d}{( b_{i}^{k} -
b_{j}^{k}) = 0}$ for all $i,j$.

But these expressions can be written even more concisely by means of scalar
products of vectors in the exterior algebra, i.e., the same technique which we
mentioned above in connection with $\tau$-critical graphs.  Towards this end,
for a vector $x = (x^{1},\ldots, x^{d}),$ we define $2^{d}$-dimensional vectors
\begin{multline*}
x^{*}  = \bigl(x^{*}(K)  \mid K \subseteq \left\{ 1,\ldots,d \right\} \bigr),  \ 
x^{\#} = \bigl(x^{\#}(K) \mid K \subseteq \left\{ 1,\ldots,d \right\} \bigr)
\\
\text{by} \quad
x^{*}(K) = \prod_{i\epsilon K}^{}x^{i} 
\quad\text{and}\quad
x^{\#}(K) = \prod_{i \notin K}^{}{- x}^{i}.
\end{multline*}
The above expressions can be then written as
\begin{equation*}
\prod_{k = 1}^{d}{(a_{i}^{k} - b_{j}^{k})} 
   = \sum \Bigl(\prod_{k \in K}^{}a_{i}^{k} \cdot \prod_{k \notin K}^{}{- b}_{j}^{k} \big\vert K \subseteq \left\{1,\ldots,d \right\} \Bigr) 
  = \sum_{K}^{}a_{i}^{*}(K) \cdot b_{j}^{\#}(K) = a_{i}^{*} \cdot b_{j}^{\#}.
\end{equation*}
Thus $a_{i}^{*} \cdot b_{j}^{\#} \neq 0$ iff $i = j.$ Similarly we
have $b_{i}^{*} \cdot a_{j}^{\#} \neq 0$ iff $i = j$ while
$a_{i}^{*} \cdot {\ a}_{j}^{\#} = b_{i}^{*} \cdot b_{j}^{\#} = 0$
for all $i,\ j.$

It follows then that the set of $2^{t}\ $vectors
\smash{$a_{i}^{*},b_{j}^{*}, \ 1 \leq i,j \leq 2^{t}$} is linearly
independent in the vector space of dimension~$2^{d}$ and thus
$t \leq d$.

Again, no other (say combinatorial) proof is known.

\section{Perfect Graphs and Computational Complexity}
\label{sec:GN7}

\begin{discussedreferences}
L. Lovász. A characterization of perfect graphs. \emph{J. Comb. Theory}
13:95--98, 1972.

L. Lovász. Normal hypergraphs and the perfect graph conjecture.
\emph{Discrete Math.} 2:253--267, 1972.
\end{discussedreferences}

\noindent
This section addresses a particular class of graphs that is tightly
connected with four important parameters. For a graph $G = (V, E)$ with
vertex set $V$ and edge set~$E$, a \emph{stable set} (also called
independent set) is a set of vertices such that no two vertices are
adjacent. The largest size of a stable set of vertices is denoted by
$\alpha(G)$ and called \emph{stability number}. Similarly, the largest size of
a \emph{clique} (mutually adjacent vertices) is denoted by $\omega(G)$ and
called \emph{clique number}, the \emph{chromatic number} $\chi(G)$ is the
smallest number of stable sets (each stable set is a color class)
covering all vertices of~$G$, and the \emph{clique covering number} $\bar{\chi}(G)$
is the smallest number of cliques covering all vertices of $G$.

If the vertices of a graph are colored so that no two adjacent vertices have
the same color then, obviously, the smallest number $\chi(G)$ of colors of such
a coloring must be at least as large as the largest number $\omega(G)$ of
mutually adjacent vertices, i.e., $\omega(G) \leq \chi(G)$. And similarly, the
stability number $\alpha(G)$ cannot be larger than the smallest number
$\bar{\chi}(G)$ of cliques covering all vertices of a graph $G$, i.e.,
$\alpha(G) \leq \bar{\chi}(G)$.

In the beginning of the 1960s Claude Berge, see~\cite{Berge1960,Berge1963},
called a graph $G$ \emph{perfect} if $\omega(H) = \chi(H)$ holds for all
induced subgraphs $H$ of $G$. In the \emph{complement} $\bar{G}$ of $G$, two
vertices are connected by an edge if and only if they are not connected in $G$,
and thus, $\alpha(G) = \omega(\bar{G})$ and $\chi(G) = \bar{\chi}(\bar{G})$.
Berge conjectured:

\begin{noheadingthm}
  A graph $G$ is perfect if and only if its complement $\bar{G}$ is perfect.
\end{noheadingthm}

This conjecture (called \emph{weak perfect graph conjecture}) started a massive
search for classes of perfect graphs. Examples are, for instance, bipartite
graphs and their line graphs, interval graphs, parity graphs, and comparability
graphs; Schrijver~\cite{Schrijver2003} describes many of these graphs in detail
in Chapter~66, Hougardy~\cite{Hougardy2006} gives a survey of these graphs and
provides a list of 120 classes. More importantly, intensive attempts to solve
the conjecture began. Fulkerson introduced pluperfect graphs
in~\cite{Fulkerson1970} and, developing in~\cite{Fulkerson1972} the
antiblocking theory for this purpose, he came very close to its solution -- as
he outlines in~\cite{Fulkerson1973}. Just a lemma (later called
\emph{replication lemma}) was missing. Lovász~\cite{Lovasz1972a} solved the
conjecture by proving the replication lemma, pointing out, though, that the
more difficult step was done first by Fulkerson. In a subsequent paper,
Lovász~\cite{Lovasz1972b} provided a new characterization of perfect graphs as
follows:

\begin{nntheorem}
A graph $G = (V, E)$ is perfect if the following holds:
$\omega(H) \alpha(H) \geq |V(H)|$ for all induced subgraphs $H = (V(H),
E(H))$ of $G$.
\end{nntheorem}

This Theorem immediately implies the weak perfect graph conjecture since
the condition given in it is invariant under taking graph
complementation. The perfect graph theorem is also a generalization of
the well-known theorems of König on bipartite matching and Dilworth on
partially ordered sets. It generated particular interest in the
characterization of conditions under which the Duality Theorem of linear
programming holds in integer variables and initiated related
investigations in polyhedral combinatorics.

Due to its importance and elegance, the Lovász's article~\cite{Lovasz1972a} was
reprinted in the collection \emph{Classic Papers in Combinatorics}, edited by
I.~Gessel and G.\,C.~Rota~\cite{GesselRota1987}.

The beginning of the 1970s was a particularly productive time period for
László Lovász. He was solving one open problem after the other. These
years firmly established his international position as the world
foremost researcher in graph theory and combinatorics.

As in many other cases, Lovász was not just looking for a proof of the
weak perfect graph conjecture, he looked for a more general mathematical
setting for which it is possible to prove farther reaching results that
imply the conjecture. In~\cite{Lovasz1972b} Lovász considered a
hypergraph approach. We sketch the construction.

Recall that a hypergraph $H$ is a non-empty finite collection of finite
sets called edges; the elements of the edges are the vertices of $H$. The
\emph{chromatic index} of a hypergraph $H$ is the least number of colors
with which the edges can be colored so that edges with the same color
are disjoint. The number of edges containing a given vertex is called
the \emph{degree} of the vertex. The largest degree of a vertex of $H$ is
called the \emph{degree of $H$}.

Clearly, the degree of $H$ is a lower bound on the chromatic index of $H$.
Lovász called a hypergraph $H$ \emph{normal} if the degree and the
chromatic index are the same for every partial hypergraph of $H$. Let us
call a set $T$ of vertices a \emph{transversal} (or hitting set) if $T$
meets every edge of $H$ and denote its minimum cardinality by $\tau(H)$. (We
just point out that $\tau(H)$ is the hypergraph generalization of $\tau(G)$ for
graphs discussed in Section~\ref{sec:GN6}.) If we denote by $\nu(H)$ the maximum number
of edges of $H$ that are pairwise disjoint, then we obviously have $\nu(H) \leq
\tau(H)$. Lovász called a hypergraph $H$ $\tau$\emph{-normal} if this inequality
holds with equality for all partial hypergraphs of $H$. He also introduced
procedures to associate with every hypergraph $H$ its edge graph~$G(H)$ and
with every graph $G$ a hypergraph $H(G)$ and proved the following:

\begin{nntheorem}
A hypergraph $H$ is normal if and only if its edge graph $G(H)$ is
perfect; $G$ is perfect if and only if $H(G)$ is normal; $H$ is $\tau$-normal if
and only if $\bar{G}(\bar{H})$ is perfect; $\bar{G}$ is perfect if and only if $H(G)$ is
$\tau$-normal.
\end{nntheorem}

\begin{nncorollary}
A hypergraph is normal if and only if it is $\tau$-normal.
\end{nncorollary}

This hypergraph generalization immediately implies the weak perfect
graph conjecture.

A side remark: In Section~\ref{sec:GN4} we mentioned the Erdős--Faber--Lovász
conjecture. This appears in this context in the following two equivalent forms:
(1)~The chromatic index of hypergraphs consisting of $n$ edges such that each
edge contains $n$ vertices and any two edges have exactly one vertex in common
is $n$. (2)~For graphs $G$ consisting of $n$ cliques of size $n$ so that two of
these cliques have one vertex in common, $\omega(G)$ equals $\chi(G)$. As
indicated before the conjecture is true for large $n$,
see~\cite{KangKellyKuehnMethukuOsthus2023}.

Berge~\cite{Berge1963} also conjectured --~later called \emph{strong
perfect graph conjecture}~-- that a graph is perfect if and only if it
does neither contain an odd cycle nor the complement of an odd cycle as
an induced subgraph. After a long sequence of contributions of many
researchers, this conjecture was finally solved in 2006 by Chudnovsky,
Robertson, Seymour, and Thomas~\cite{ChudnovskyRobertsonSeymourThomas2006}.

During the early 1970s computational complexity theory took off, see
Wigderson's book~\cite{Wigderson2019} for an up-to-date survey. The classes of
decision problems that can be solved in polynomial time, denoted by \cP,
and those solvable in nondeterministic polynomial time, denoted by \cNP,
were introduced. S.~Cook~\cite{Cook1971} and L.\,A.~Levin~\cite{Levin1973}
independently showed the existence of \emph{\cNP-complete} problems,
which are decision problems in \cNP with the property that, if they can
be solved with a polynomial time algorithm, then $\cP = \cNP$. Whether
\cP is equal to \cNP is one of the great open problems in
mathematics and computer science.

Optimization problems can be phrased as decision problems by asking whether,
for a given value~$t$, there exists a feasible solution with value at least (or
at most)~$t$. If the decision problem associated this way to an optimization
problem is \cNP-complete, the optimization problem is called
\emph{\cNP-hard}. For example, if a graph \mbox{$G = (V, E)$} with rational
weights~$w_{v}$ for every vertex $v \in V$, is given and one wants to find a
stable set $S$ in $V$ such that the sum of the weights of the vertices in $S$
is as large as possible, we have a typical combinatorial optimization problem.
The associated decision problem asks if there is a stable set whose value is at
least $t$. If this decision problem can be solved in polynomial time, the
stable set problem can also be solved in polynomial time by binary search. And
vice versa, a polynomial time algorithm for the (weighted) stable set problem
would prove that $\cP = \cNP$.

Karp~\cite{Karp1972} showed that many graph-theoretical problems, such as
computing the value of the four parameters $\alpha(G)$, $\omega(G)$, $\chi(G)$,
and $\bar{\chi}(G)$, introduced above, are \cNP-hard for general graphs
$G$. The immediate question came up: Is that also true for perfect graphs, or
can their special structure be exploited to design polynomial time algorithms?
This challenge triggered significant developments that we outline later.

Another side remark: Lovász was one of many contributors to one of the most
astonishing results in complexity theory, the PCP Theorem. This theorem is the
highlight of a long sequence of research on interactive proofs and
probabilistically checkable proofs. It states that every decision problem in
\cNP has probabilistically checkable proofs of constant query complexity
using only a logarithmic number of random bits. Nine persons (including Lovász)
received the Gödel Prize~2002 ``\emph{for the PCP theorem and its applications
to hardness of approximation}''. A consequence of the PCP Theorem is, for
instance, that many well-known optimization problems, including the stable set
problem mentioned above and the shortest vector problem for lattices to be
introduced subsequently, cannot be approximated efficiently unless $\cP =
\cNP$.

\section{The Shannon Capacity of a Graph and Orthogonal Representations}
\label{sec:GN8}

\begin{discussedreferences}
  L. Lovász. On the Shannon capacity of graphs. \emph{IEEE Trans. Inform.
  Theory} 25:1--7, 1979.

L. Lovász. Graphs and geometry. \emph{Amer. Math. Soc.} 2019.
\end{discussedreferences}

\noindent
Suppose the vertices of a graph $G$ represent letters of an alphabet and
the edges~$uv$ of $G$ indicate that the two letters of the alphabet
represented by $u$ and $v$ can be confused, e.g., when transmitted over a
noisy communication channel. It is obvious that the largest number of
one-letter messages that can be sent without danger of confusion is the
largest number of vertices mutually not adjacent, i.e., the stability
number $\alpha(G)$. Two $k$-letter words are confusable if their $i$-th letters, $1
\leq i \leq k$, are confusable or equal.

Let $G^{k}$ denote the $k$-th Cartesian product of $G$. Words
with $k$-letters can be transmitted without danger of confusion if they
are unequal and inconfusable in at least one letter. This implies that
$\alpha(G^{k})$ is the maximum number of inconfusable $k$-letter
words. Forming $k$-letter words from a stable set of size $\alpha(G)$, one can
easily construct $\alpha(G)^{k}$ inconfusable words. This proves
that $\alpha(G)^{k} \leq \alpha(G^{k})$.

\begin{figure}[t]
\centering
\includegraphics[width=4cm]{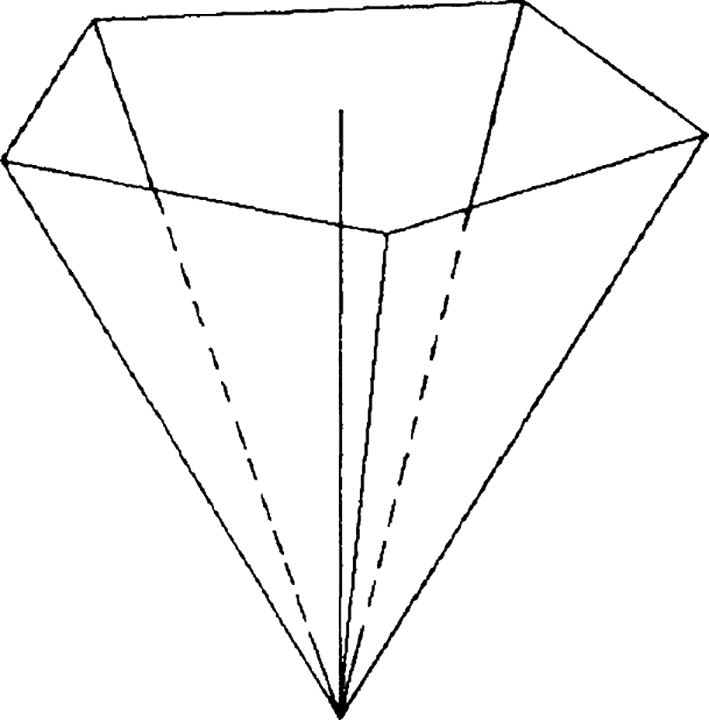}

\caption{Orthonormal representation of the 5-cycle in $\mathbb{R}^3$}
\label{fig:GN7}
\end{figure}

Shannon~\cite{Shannon1956} introduced the number
\[
  \Theta(G)
  = \sup_k \sqrt[k]{\alpha(G^k)} 
  = \lim_{k \rightarrow \infty} \sqrt[k]{\alpha(G^k)},
\]
where the second equation follows from $\alpha(G^{k+l}) \geq \alpha(G^{k})
\alpha(G^{l})$. $\Theta(G)$, today called the \emph{Shannon capacity} of $G$,
is a measure of the information that can be transmitted across a noisy
communication channel. Shannon proved that $\Theta(G) = \alpha(G)$ for graphs
which can be covered by $\alpha(G)$ cliques. Perfect graphs have this property
and thus belong to this class. How can one determine $\Theta(G)$ in other
cases? Lovász, see \cite{Lovasz1979a}, invented an ingenious upper bound on the
Shannon capacity as follows:

Let $G = (V, E)$ be a graph. An \emph{orthonormal representation} of $G$ is a
sequence $(u_i \mid i \in  V)$ of $|V|$ vectors $u_i \in
\mathbb{R}^N$, where $N$ is some positive integer, such that $ \|u_i\| = 1$ for
all $i \in V$ and $u_i^T u_j = 0$ for all pairs $i,j$ of nonadjacent vertices.
Trivially, every graph has an orthonormal representation (just take all the
vectors $u_i$ mutually orthogonal in $\mathbb{R}^V$). Figure~\ref{fig:GN7}
shows a less trivial orthonormal representation of the pentagon $C_5$ in
$\mathbb{R}^3$. It is constructed as follows.  Consider an umbrella with five
ribs of unit length (representing the nodes of $C_5$) and open it in such a way
that nonadjacent ribs are orthogonal. Clearly, this can be achieved in
$\mathbb{R}^3$ and gives an orthonormal representation of the pentagon.  The
central handle (of unit length) is also shown.

Where $(u_i \mid i \in V), u_i \in \mathbb{R}^N$, ranges over all orthonormal
representations of $G$ and $c \in \mathbb{R}^N$ over all vectors of unit
length, let
\[
  \vartheta (G, w) \coloneqq \min_{\{c, (u_i)\}} \max_{i \in V} \frac{w_i}{(c^T u_i)^2}.
\]
The quotient has to be interpreted as follows. If $w_i = 0$ then we take
$w_i/(c^T u_i)^2 = 0$ even if $c^T u_i = 0$. If $w_i > 0$ but $c^T u_i = 0$
then we take $w_i/(c^T u_i)^2 = +\infty$.

Lovász proved that, if the vertex weights $w_{i}$ above are all equal to~1 and
$G$ is the pentagon graph $C_{5}$, i.e., the 5-cycle, then the value of
$\vartheta(G,w)$ is $\sqrt{5}$ and equal to the Shannon capacity $\Theta(C_5)$ of $C_{5}$.

This looks like a tiny achievement, but at present, this is the only known
Shannon capacity of a non-perfect graph. In fact, the complexity of determining
the Shannon capacity of a general graph is today still open. Much more
important, Lovász provided several different characterizations of the function
$\vartheta$ (called the \emph{Lovász $\vartheta$\nobreakdash-function}) that
became, as we show later, important ingredients for proving that the four graph
parameters $\alpha(G)$, $\omega(G)$, $\chi(G)$, and $\bar{\chi}(G)$ can be
computed in polynomial time for perfect graphs $G$.

In his recent book~\cite{Lovasz2019}, Lovász investigated the
representation of graphs as geometric objects in great depth. His main
message is that such representations are not merely a way to visualize
graphs, but important mathematical tools. The range of applications is
wide. We mention three examples: rigidity of frameworks and mobility of
mechanisms in engineering, learning theory in computer science, the
Ising and Fortuin--Kasteleyn model, and conformal invariance in
statistical physics. Orthogonal representations of graphs are treated in
Chapters~10 to 12. Lovász shows that orthogonal representations are, in
addition to the stability and chromatic number, related to several
fundamental properties of graphs such as connectivity and tree-width.
Among many other aspects, he also discusses a quantum version of the
Shannon capacity problem, as well as two further interesting
applications of orthogonal representations to the theory of hidden
variables and in the construction of strangely entangled states. These
are exciting topics in quantum physics that we cannot cover here.

\newpage

\section{The Ellipsoid Method}
\label{sec:GN9}

\begin{discussedreferences}
  P. Gács, L. Lovász. Khachiyan's algorithm for linear programming. \emph{Math.
  Prog. Study} 14: 61--68, 1981.
\end{discussedreferences}

\noindent
One of the major open complexity problems in the 1970s was the question
whether linear programs (LPs) can be solved in polynomial time. The
simplex algorithm did (and still does) work well in practice, but for
all known variants of this algorithm, there exist sequences of
LP-instances for which the running time is exponential. In~1979
Khachiyan indicated in~\cite{Khachiyan1979} how the ellipsoid method, an
algorithm devised for nonlinear nondifferentiable optimization based on
work of Shor and Yudin and Nemirovski\u{\i}, can be modified to check the
feasibility of a system of linear inequalities in polynomial time.
Employing binary search or a sliding objective function technique, this
implies that linear programs are solvable in polynomial time. Linear
programs arise almost everywhere in industry, and their fast solution is
of economic importance. Thus, Khachiyan's achievement received
significant attention in the nonscientific media; it even made it on the
front page of the New York Times on November~7, 1979. Most of these
statements, though, were exaggerations or misinterpretations.

We sketch the method. Let $P$ be polyhedron defined by a system of linear
inequalities $Ax\leq b$. We assume that $P$ is full-dimensional or empty; and
for simplifying the exposition, we also assume that $P$ is bounded, i.e.,
a polytope. The ellipsoid method utilizes the following facts. Given
$Ax\leq b$ with rational coefficients, then numbers~$r$ and $R$ can be computed in
time polynomial in the encoding length of $A$ and $b$ with the following
properties. If $P$ is nonempty, the ball~$B$ of radius~$R$ around the origin
contains $P$, and $P$ contains a ball~$S$ of radius~$r$.

\begin{figure}
\centering
\includegraphics[width=5cm]{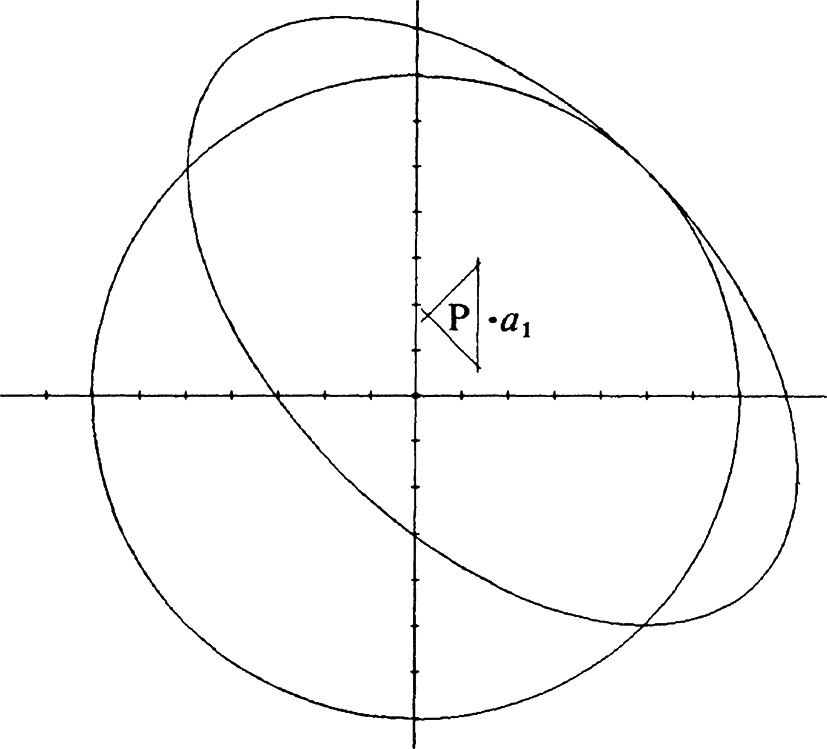}

\caption{The first step of the ellipsoid method}
\label{fig:GN8}
\end{figure}

The basic ellipsoid method begins with the ball $B$ and center $a_{0} = 0$ as
initial ellipsoid $E_{0}$. In a general step it checks whether the center
$a_{k}$ of the current ellipsoid~$E_{k}$, $0 \leq k$, is contained in $P$. If
this is the case, a point in $P$ is found and $Ax\leq b$ is feasible. If not,
there must be an inequality in the system $Ax\leq b$ that is violated by
$a_{k}$. Using this inequality, a new ellipsoid $E_{k+1}$ is computed that
contains $P$ and has a volume that is -- by a constant shrinking rate --
smaller than the volume of the previous ellipsoid $E_{k}$
(cf.~Fig.~\ref{fig:GN8}). This way a sequence of points $a_{k}$ and shrinking
ellipsoids~$E_{k}$ is created. Using variants of the formulas for determining
the Löwner--John-ellipsoid of a convex body, one can prove that the volume
shrinking rate satisfies $\vol(E_{k+1})/\vol(E_{k}) < e^{-1/(2n)} < 1$ and that
the ellipsoid method either discovers a point in $P$ or, after a number $N$ of
steps that is polynomial in the encoding length of~$A$ and $b$, the ellipsoid
$E_{N}$ has a volume that is smaller than that of the small ball $S$. This can
only happen in case $P$ is empty. All computations carried out can be made with
rational numbers of polynomial size in such a way that nonemptiness of $P$ is
certified by finding a feasible solution or the emptiness of $P$ is
guaranteed by the mentioned volume argument, see~\cite{GLS1988} for details.

This method was a total surprise for the linear programming community. A
polynomial time termination proof employing shrinking volumes, the combination
of geometric and number theoretic ``tricks'' (e.g., making a low-dimensional
polyhedron full-dimensional, reduction to the bounded case, careful rounding of
the real numbers that appear in the update-formulas, and various necessary
estimation processes) puzzled the LP-specialists. The brief article by
Khachiyan (four pages), written in Russian, needed interpretation. One of the
first papers explaining the approach and adding missing details was a preprint
by Gács and Lovász~\cite{GacsLovasz1981}. It appeared in the fall of 1979 (and
was published in~1981). This paper made Khachiyan's important contribution
accessible to a wide audience and had a significant bearing on the boom of
follow-up research on the ellipsoid method.

The ellipsoid method, though provably a polynomial time algorithm,
performs poorly in practice. Its appearance, however, sparked successful
research efforts that led to new LP-algorithms, based on various ideas
from nonlinear programming, often also influenced by differential and
other types of geometry, that are theoretically and practically fast.
They run under the names interior point or barrier methods. New
implementations of the simplex algorithm improved its performance
significantly as well. The ellipsoid method, on the other hand, turned
out to have fundamental theoretical power as an elegant and versatile
tool to prove the polynomial time solvability of many geometric and
combinatorial optimization problems. The next chapter has details.

\section{Oracle-Polynomial Time Algorithms and Convex Bodies}
\label{sec:GN10}

\begin{discussedreferences}

M. Grötschel, L. Lovász, A. Schrijver. The ellipsoid method and its
consequences in combinatorial optimization, \emph{Combinatorica}
1:169--197, 1981.

M. Grötschel, L. Lovász, A. Schrijver. \emph{Geometric Algorithms and
Combinatorial Optimization}, Springer, 1988.
\end{discussedreferences}

\noindent
In a general step of the ellipsoid method, one has to verify that the center of
the current ellipsoid is in the polyhedron $P=\{x \in \mathbb{R}_{n} \mid
Ax\leq b\}$. This is usually done by substituting the center into the given
inequality system $Ax\leq b$. A reasonable idea is to replace this substitution
by an algorithm that checks feasibility and provides a violated inequality in
case the center is not in $P$. Two cases, relevant in real-world applications,
where this generalization might be helpful come immediately into mind.

\begin{figure}[b]
\centering
\includegraphics[width=\linewidth]{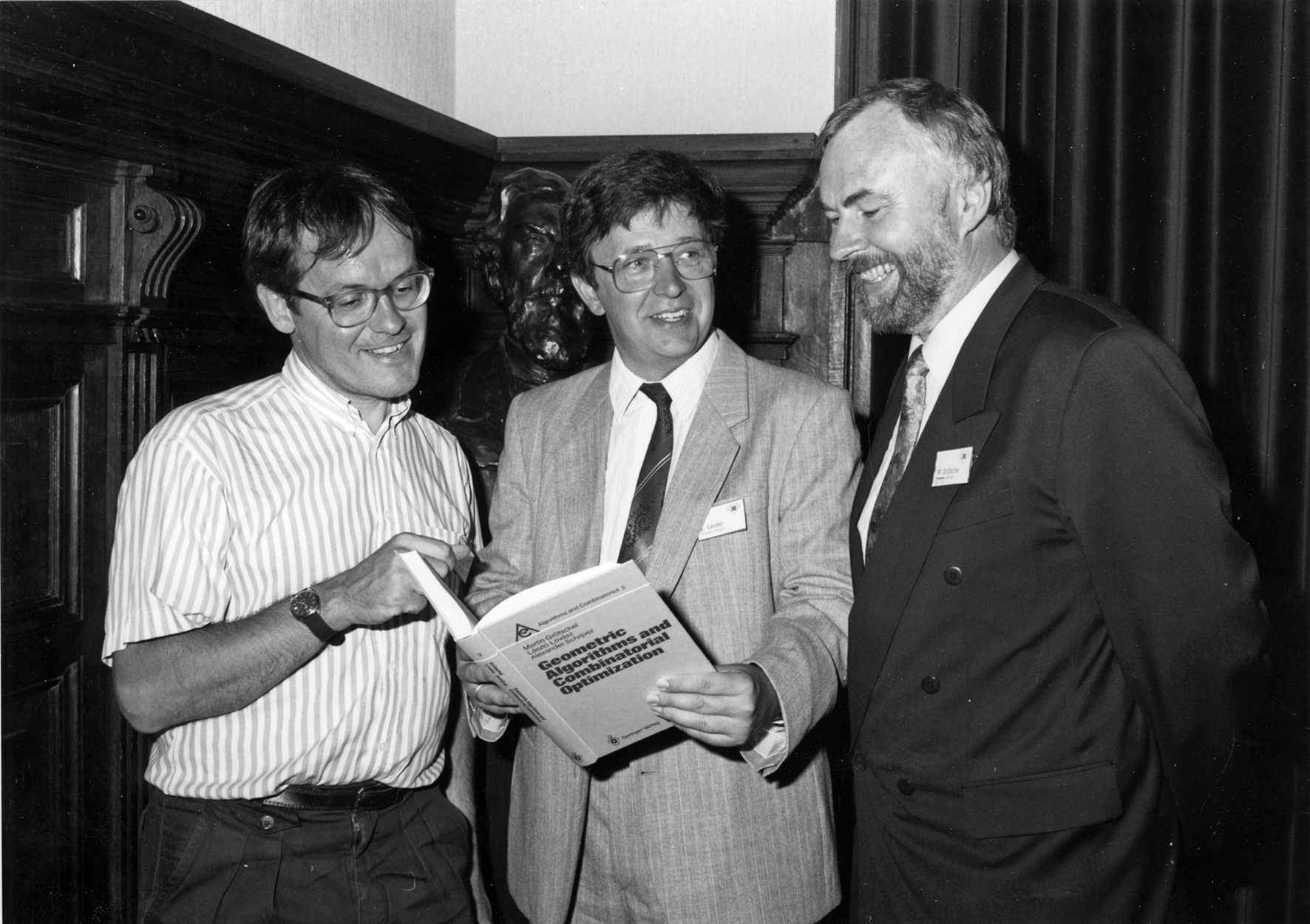}

\caption{A. Schrijver, L. Lovász, M. Grötschel at the International Symposium
on Mathematical Programming in Amsterdam, 1991 (Photo: Nationaal
Foto-Persbureau B.\,V.)}
\label{fig:9}
\end{figure}

The first case is the traditional transformation of combinatorial optimization
problems into linear programs. The idea is, for a given combinatorial
optimization problem, to define the convex hull of all incidence vectors of
feasible solutions and to try to find a linear system describing this polytope,
at least partially. The number of facets of such polytopes is often
exponentially large in the encoding length of the combinatorial problem. This
holds for \cNP-hard problems and even for some problems solvable in
polynomial time. One such instance is the matching problem. This is implied by
the result of Rothvoss~\cite{Rothvoss2017} that the matching problem has
``exponential extension complexity''. Substituting the ellipsoid center into a
linear system of exponential size makes the running time of ellipsoid algorithm
exponential. Can one replace the substitution by a polynomial time algorithm?

The second case are convex sets and is even more demanding. Convex sets are
intersections of potentially infinitely many halfspaces. Can one optimize over
exponentially many linear inequalities in polynomial time?

The roots of this research program were laid by Grötschel, Lovász, and
Schrijver in~\cite{GLS1981} and were fully worked out in~\cite{GLS1988}.  The
results were the starting point of what Gritzmann and
Klee~\cite{GritzmannKlee1997} called an algorithmic theory of convex bodies, or
briefly, \emph{computational convexity}. We outline important steps of this
approach.

Suppose now that we have some convex set $K \subseteq \mathbb{R}_{n}$ and we
want to obtain information about properties of $K$. Let us formulate three
questions that are typical in this context:

\begin{SOPT}
Given a vector $c \in \mathbb{R}^n$, find a vector $y \in K$ that maximizes
$c^T x$ on $K$, or assert that K is empty.
\end{SOPT}

\begin{SSEP}
Given a vector $y \in \mathbb{R}^n$, decide whether $y \in K$, and if not, find
a hyperplane that separates $y$ from $K$; more exactly, find a vector $c \in
\mathbb{R}^n$ such that $c^T y > \max\{c^T x \mid x \in K\}$.
\end{SSEP}

\begin{SMEM}
Given a vector $y \in \mathbb{R}^n$, decide whether $y \in K$.
\end{SMEM}

It is clear that the strong membership problem can be solved if either
the strong optimization or the strong separation problem can be solved.
What about the other way around? And what do we have to assume about~$K$,
what is the input length of~$K$, and how do we estimate running times?
Before addressing these issues, we observe that, if we allow arbitrary
convex sets~$K$, the unique solution of an optimization problem over~$K$ may
have irrational coordinates. To deal with such issues we have to allow
margins and to accept approximate solutions. Let us define, for the
Euclidean norm and a rational number $\epsilon > 0$,
\[
 S (K, \epsilon)  \coloneqq \bigl\{x \in \mathbb{R}^n \mid \|x-y\| \leq \epsilon \text{ for some } y \in K\bigr\} , \
 S (K, -\epsilon) \coloneqq \bigl\{x \in K \mid S(x,\epsilon) \subseteq K\bigr\}.
\]
Points in $S(K, \epsilon)$ can be viewed as ``almost in $K$'', while points in $S(K,
-\epsilon)$ as ``deep in $K$''. The exactness requirements of the strong problems
above can be softened as follows:

\begin{WOPT}
Given a vector $c \in \mathbb{Q}^n$ and a rational number $\epsilon > 0$,
either
\begin{compactenum}[\rm(i)]
\item find a vector $y \in \mathbb{Q}^n$ such that $y \in S(K, \epsilon)$ and
  $c^Tx \leq c^Ty + \epsilon$ for all $x \in S(K, -\epsilon)$ (i.e., $y$ is
  almost in $K$ and almost maximizes $c^Tx$ over the points deep in $K$), or
\item assert that $S(K, -\epsilon)$ is empty.
\end{compactenum}
\end{WOPT}

\begin{WSEP}
Given a vector $y \in \mathbb{Q}^n$ and a rational number $\delta > 0$, either
\begin{compactenum}[\rm(i)]
\item assert that $y \in S(K, \delta)$, or
\item find a vector $c \in \mathbb{Q}^n$ with $\|c\|_\infty = 1$ such that $c^T
  x \leq c^T y + \delta$ for every $x \in S(K, -\delta)$ (i.e., find an almost
  separating hyperplane).
\end{compactenum}
\end{WSEP}

\begin{WMEM}
Given a vector $y \in \mathbb{Q}^n$ and a rational number $\delta > 0$, either
\begin{compactenum}[\rm(i)]
\item assert that $y \in S (K,\delta)$, or 
\item assert that $y \notin S (K, -\delta)$.
\end{compactenum}
\end{WMEM}

We are interested in the algorithmic relations between these problems.  To do
this we make use of the \emph{oracle algorithm concept}. An \emph{oracle} is a
device that solves a certain problem for us. Its typical use is as follows. We
feed some input string to the oracle, and the oracle returns another string
specifying the solution (of which we hope that it helps solving our original
problem). We make no assumption on the way the oracle finds its solution. An
oracle algorithm is an algorithm in the usual sense whose power is enlarged by
allowing querying an oracle and using the oracle answer for determining its
next computational steps.

If a query to and an answer of the oracle are counted as one step each, we can
determine the running time of an oracle algorithm in the usual way. The output
of the oracle may, however, be huge so that reading it may take exponential
time. Since our aim is to design polynomial time algorithms, we require that
for every oracle we have a polynomial $q$, such that for every query of
encoding length at most $l$, the answer of the oracle has length at most
$q(l)$. Under this assumption we say that an oracle algorithm has
\emph{oracle-polynomial running time} if its usual running time plus the
running time of the interaction with the oracle is bounded by a polynomial in
the input length of the original problem. A consequence of this set-up is that,
if an oracle can be realized by a polynomial time algorithm on a real
computational device, an oracle-polynomial algorithm is in fact a polynomial
time algorithm in the usual sense.

For ease of exposition, we restrict ourselves to considering convex bodies~$K$
only. A convex set~$K \subseteq \mathbb{R}^{n}$ that is compact and has
dimension~$n$ is called \emph{convex body}. To perform computations, we have to
assume that the convex body~$K$ is given by a mathematical description.  Let us
briefly call it $\Name(K)$. Then the encoding length of~$K$ is defined as the
dimension~$n$ plus the encoding length of~$\Name(K)$. To determine the
algorithmic relations between the problems above, we assume that a convex body
is given by an oracle for the solution of one of the problems and we
investigate whether any of the other problems can be solved employing the
oracle. The running times are measured as usual in the size of the input. This
is, in the cases described here, the encoding length of~$K$ (as defined above)
to which we have to add, if they appear in the problem statement, the
following: the encoding lengths of the parameters~$\epsilon$ and $\delta$, the
encoding lengths of the objective function~$c$ and the vector~$y$, and moreover
the encoding lengths of the additional data (the radii~$r$ and $R$, and the
center $a_{0}$ of a ball) appearing in the statements of the theorems. The
following was proved in~\cite{GLS1988}:

\begin{nntheorem}
\begin{compactenum}[\rm(a)]
\item There exists an-oracle polynomial time algorithm that solves the weak
  membership problem for every convex body $K$ in $\mathbb{R}^{n}$ given by a
  weak optimization or a weak separation oracle.

\item There exists an oracle-polynomial time algorithm that solves the weak
  separation problem for every convex body $K$ in $\mathbb{R}^{n}$ given by a
  weak optimization oracle.

\item There exists an oracle-polynomial time algorithm that solves the weak
  optimization problem for every convex body $K$ in $\mathbb{R}^{n}$ given by a
  weak separation algorithm, provided a radius $R>0$ of a ball around the
  origin containing $K$ is given as well.

\item There exists an oracle-polynomial time algorithm that solves the weak
  optimization problem for every convex body $K$ in $\mathbb{R}^{n}$ given by a
  weak membership algorithm, provided the following data are given as well: a
  vector $a_{0}$ and a radius $r>0$ such that
  $S(a_{0},r) \subseteq K$, and a radius $R>0$ with
  $K\subseteq S(0,R)$.
\end{compactenum}
\end{nntheorem}

This theorem establishes the oracle-polynomial time equivalence of WOPT, WSEP,
and WMEM under mild additional assumptions. Moreover, the oracle-polynomial
time equivalence of the strong versions SOPT, SSEP, and SMEM of these problems
can be derived from the results above (assuming, of course, that $K$ is given
such that exact answers are possible). One can prove on the other hand that, if
we would drop one of the additional requirements in the theorem such as the
knowledge of radii~$r$ or $R$ or the vector~$a_{0}$, it is impossible to derive
oracle-polynomial time algorithms.

A consequence of the last result, see \cite{GLS1981} and \cite{GLS1988}, is the
polynomial time solvability of convex function minimization -- in the following
weak sense:

\begin{nntheorem}
There exists an oracle-polynomial time algorithm that solves the
following problem:
\\
\emph{Input:} A convex body $K$ given by a weak membership oracle,
a rational number $\epsilon > 0$, radii~$r$, $R>0$, a vector
$a_{0}$ such that $S(a_{0}, r) \subseteq K \subseteq S(0, R)$, and a
convex function $f : \mathbb{R}^{n} \rightarrow \mathbb{R}$ given by an
oracle that, for every $x \in \mathbb{Q}^n$ and $\delta > 0$, returns a
rational number $t$ such that $|f(x) - t| < \delta$.
\\
\emph{Output:} A vector $y \in S(K,\epsilon)$ such that $f(y) < f(x) +
\epsilon$ for all $x \in S(K,-\epsilon)$.
\end{nntheorem}

This is the first polynomial time solvability result for convex
minimization.

\section{Polyhedra, Low Dimensionality, and the LLL Algorithm}
\label{sec:GN11}

\begin{discussedreferences}

M. Grötschel, L. Lovász, A. Schrijver. The ellipsoid
method and its consequences in combinatorial optimization.
\emph{Combinatorica}, 1:169--197, 1981.

A.\,K. Lenstra, H.\,W. Lenstra, L. Lovász. Factoring Polynomials with
rational coefficients. \emph{Mathematische Annalen} 261(4):515--534,
1982.
\end{discussedreferences}

\noindent
The Abel prize citation states (correctly, of course): ``\emph{The LLL
algorithm is only one among many of Lovász's visionary contributions}''.  It
may be surprising to learn that its invention was triggered by a technical
problem arising in the analysis of the ellipsoid method. We explain its origin
and usefulness in this context.

Since square roots appear in the update formulas defining the ellipsoid method,
computing with irrational numbers is unavoidable. Careful rounding is necessary
to reach the desired approximation of an optimal value or solution. In various
applications exact solutions can in fact be obtained by appropriate rounding.
In integer programming, e.g., the solution vectors are required to have
integral entries, and if the objective function is integral, the optimal
value~$v^*$ is integral as well. If one can tune the ellipsoid method so that
it guarantees to find an approximation $v$ of the optimal value~$v^*$ such that
$|v-v^*| < 1/2$, then one can simply round $v$ to the next integer to find the
true optimum value. Such considerations are the key to pass from ``weak
solutions'' to ``strong solutions'', i.e., derive exact from approximate
results. This straightforward rounding unfortunately is often not sufficient.

We sketch the case of optimizing a linear objective function over a
polytope~\mbox{$P \subseteq \mathbb{R}^{n}$}. We say that $P$ has
\emph{facet-complexity at most $\phi$} if there exists a system of
inequalities with rational coefficients that has solution set $P$ and such
that the encoding length of each inequality of the system is at most $\phi$.
No assumption about the number of inequalities is made. Let us define
the encoding length of~$P$ to be $n + \phi$, call such a polyhedron
\emph{well-described,} and denote it by $(P; \ n,\phi)$. One can prove that the
encoding length of each vertex of~$(P; \ n,\phi)$ is at most
$4n^{2}\phi$ and that, if $P$ is full-dimensional, $P$ contains a
ball~$B_{P}$ with radius~$2^{-7n^3\phi}$.

To illustrate the annoying ``technical problem'' that triggered the invention
of the LLL algorithm, let us consider a well-described polytope $P \subseteq
\mathbb{R}^{n}$ that is not full-dimensional; for ease of exposition, say $P$
has dimension $n-1$. The ellipsoid method would not work in this case. To get
around this problem, one needs to carefully blow~$P$ up to a polytope~$P'$ that
contains~$P$ and is full-dimensional such that running the ellipsoid method
on~$P'$ approximately delivers the desired result for $P$. This can be done but
is technically tedious and requires ugly pre- and post-processing.

Let us instead make a bold step and run the ellipsoid method on $P$ directly.
We suppose $P$ is given by a separation oracle. Since $P$ is low-dimensional it
is highly unlikely that the ellipsoid method finds a feasible solution in one
of its iterations. After a number $N$ of iterations that is polynomial in $n +
\phi$, the $N$-th ellipsoid $E_{N}$ contains~$P$ and has a volume that is
smaller than the volume of~$B_{P}$, the ball~$P$ would contain if $P$ were
full-dimensional. This is contradictory. The basic ellipsoid method, assuming a
full-dimensional polytope~$P$ is given, would conclude now that $P$ is empty.
But $E_{N}$ contains information that one may be able to employ.

Let $H = \{x \in \mathbb{R}^{n} \mid a^{T}x = \alpha\}$ be the unique hyperplane
containing~$P$. Then $a^{T}x = \alpha$ is the (up to scaling) unique equation
defining~$H$. The last ellipsoid~$E_{N}$, having such a small
volume, must obviously be very ``flat'' in the direction perpendicular to~$H$.
In other words, the symmetry hyperplane~$F$ belonging to the shortest axis of
$E_{N}$ must be very close to~$H$. Is it possible to find $a^{T}x = \alpha$ by
rounding the coefficients of the linear equation defining this symmetry
hyperplane~$F$? A positive answer would be an elegant way to avoid the blow-up
mentioned and the numerical problems associated with it.

The authors of~\cite{GLS1981} and~\cite{GLS1988} were at this point in the fall
of 1981 and realized that such a rounding can be done -- in principle -- using
the following classical theorem of Dirichlet~\cite{Dirichlet1842} on the
existence of a solution of a simultaneous Diophantine approximation problem.

\begin{nntheorem}
Given any real numbers $\alpha_{1}, \ldots , \alpha_{n}$ and $0 < \epsilon <
1$, there exist integers $p_{1}, \ldots , p_{n}$, and $q$ such that $1 < q <
\epsilon^{-n}$ and $|\alpha_{i} - p_{i}/q| < \epsilon/q$ for $i = 1, \ldots , n$.
\end{nntheorem}

No polynomial algorithm is known to compute such integers. And at the end of
their writing session, no progress was achieved. About three months later the
following letter from L.~Lovász arrived:

\begin{figure}[H]
\centering
\includegraphics[width=.85\linewidth]{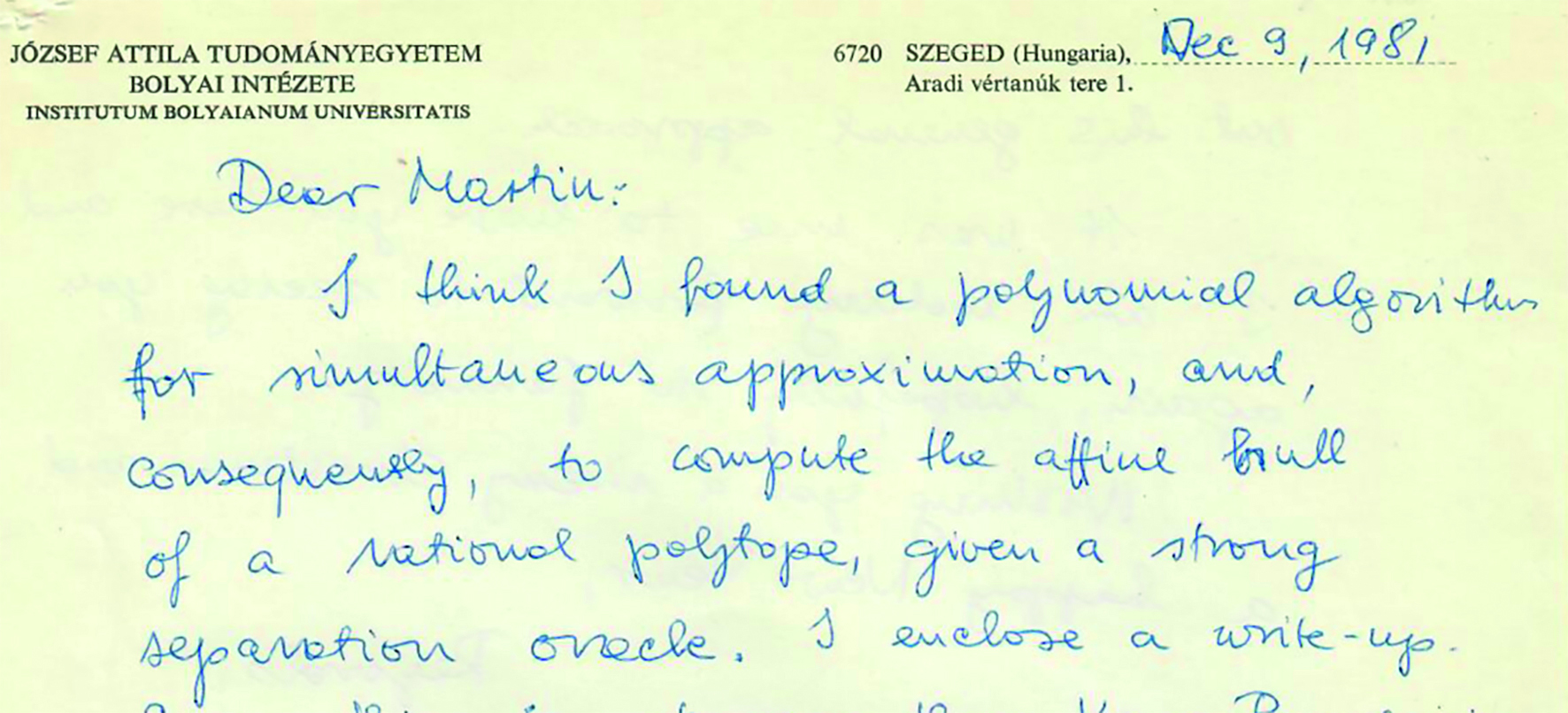}

\caption{Beginning of a letter from L. Lovász}
\label{fig:GN10}

\end{figure}

Lovász approached the approximation problem via the consideration of (integral)
lattices. If $\{b_{1}, \ldots , b_{n}\}$ is a basis of $\mathbb{R}^{n}$, then
the set $L = L(b_{1}, \ldots , b_{n})$ that is generated by taking all integral
linear combinations of the vectors $b_{i}$ is called a \emph{lattice} with
basis $\{b_{1}, \ldots , b_{n}\}$. Integral lattices have been studied in
number theory for a very long time (with contributors such as Gauss, Minkowski,
Landau, and many others). Clearly, a lattice may have different bases, and it
may be interesting to find a ``minimal basis'' $\{a_{1}, \ldots , a_{n}\}$ of
$L$, i.e., a basis such that the product of the norms of the $a_{i}$ is as
small as possible. However, this problem is \cNP-hard. Lovász introduced
the quite technical notion of a \emph{reduced basis}, which we do not explain
here, that is a weak form of a minimal basis and proved:

\begin{nntheorem}
There is a polynomial time algorithm that, for any given linearly independent
vectors $\{b_{1}, \ldots , b_{n}\}$ in $\mathbb{Q}^{n}$, finds a reduced basis
of the lattice $L(b_{1}, \ldots , b_{n})$.
\end{nntheorem}

The algorithm, called LLL algorithm, to achieve this starts with the
Gram--Schmidt orthogonalization and then performs carefully designed exchange
operations. Proving polynomiality requires not only controlling the number of
steps, but in particular, the estimation of the encoding lengths of all numbers
appearing in the course of the algorithm. A consequence of this algorithm is
the following weak form of Dirichlet's theorem.

\begin{nntheorem}
There exists a polynomial time algorithm that, given rational numbers
$\alpha_{1}, \ldots , \alpha_{n}$ and $0 < \epsilon < 1$, computes integers
$p_{1}, \ldots , p_{n}$, and $q$ such that and $1 \leq q \leq
2^{n(n+1)/4}\epsilon^{-n}$ and $|\alpha_{i}q - p_{i}| < \epsilon$ for $i = 1,
\ldots , n$.
\end{nntheorem}

This algorithm, based on computing a reduced basis, made it possible to compute
via simultaneous Diophantine approximation, the coefficients of the
equation~$a^{T}x = \alpha$ defining the hyperplane~$H$ containing the
well-described polytope $(P; \ n,\phi)$ as indicated above. By iterating this
process, the affine hull of any lower dimensional polytope can be determined in
oracle-polynomial time.

For well-described polyhedra $(P; \ n,\phi)$, the restriction to the bounded
case can also be dropped, and one can show the following:

\begin{nntheorem}
Any of the following three problems:
\begin{compactitem}[--]
\item strong separation
\item strong violation
\item strong optimization
\end{compactitem}
can be solved in oracle-polynomial time for any well-described polyhedron $(P;
\ n,\phi)$ given by an oracle for any of the other two problems.
\end{nntheorem}

For a linear program given by a system of rational linear inequalities, the
strong separation problem can be trivially solved by substituting a given
rational vector y into the inequalities, i.e., linear programs can be solved in
polynomial time.

Employing the LLL algorithm and results of András Frank and Éva
Tardos~\cite{FrankTardos1987} one can, in fact, derive a general result about
optimization problems for polyhedra and their dual problems in strongly
polynomial time. \emph{Strongly polynomial} means that the number of elementary
arithmetic operations to solve an optimization problem over a well-described
polyhedron and to solve its dual problem does not depend on the encoding length
of the objective function. More precisely, the following can be shown:

\begin{nntheorem}
There exist algorithms that, for any well-described polyhedron $(P; \ n,\phi)$
specified by a strong separation oracle, and for any given vector $c \in
\mathbb{Q}^{n}$,
\begin{compactenum}[\rm(a)]
\item solve the strong optimization problem $\max\{c^{T}x \mid x \in P\}$, and
\item find an optimum vertex solution of $\max\{c^{T}x \mid x \in P\}$ if one exists, and
\item find a basic optimum standard dual solution if one exists. 
\end{compactenum}
The number of calls on the separation oracle, and the number of elementary
arithmetic operations executed by the algorithms are bounded by a polynomial
in~$\phi$. All arithmetic operations are performed on numbers whose encoding
length is bounded by a polynomial in~$\phi$ and the encoding length of the
objective function vector~$c$.
\end{nntheorem}

An important application of this theorem is that one can turn many polynomial
time combinatorial optimization algorithms into strongly polynomial algorithms.

Summarizing: The search for an elegant proof that avoids tedious numerical
estimates was the driving force for the invention of the LLL algorithm.

\section{The LLL Algorithm and its Consequences}
\label{sec:GN12}

\begin{discussedreferences}
A.\,K. Lenstra, H.\,W. Lenstra, L. Lovász. Factoring polynomials with
rational coefficients. \emph{Mathematische Annalen} 261(4):515--534,
1982.
\end{discussedreferences}

\begin{figure}
\centering
\includegraphics[width=.85\linewidth]{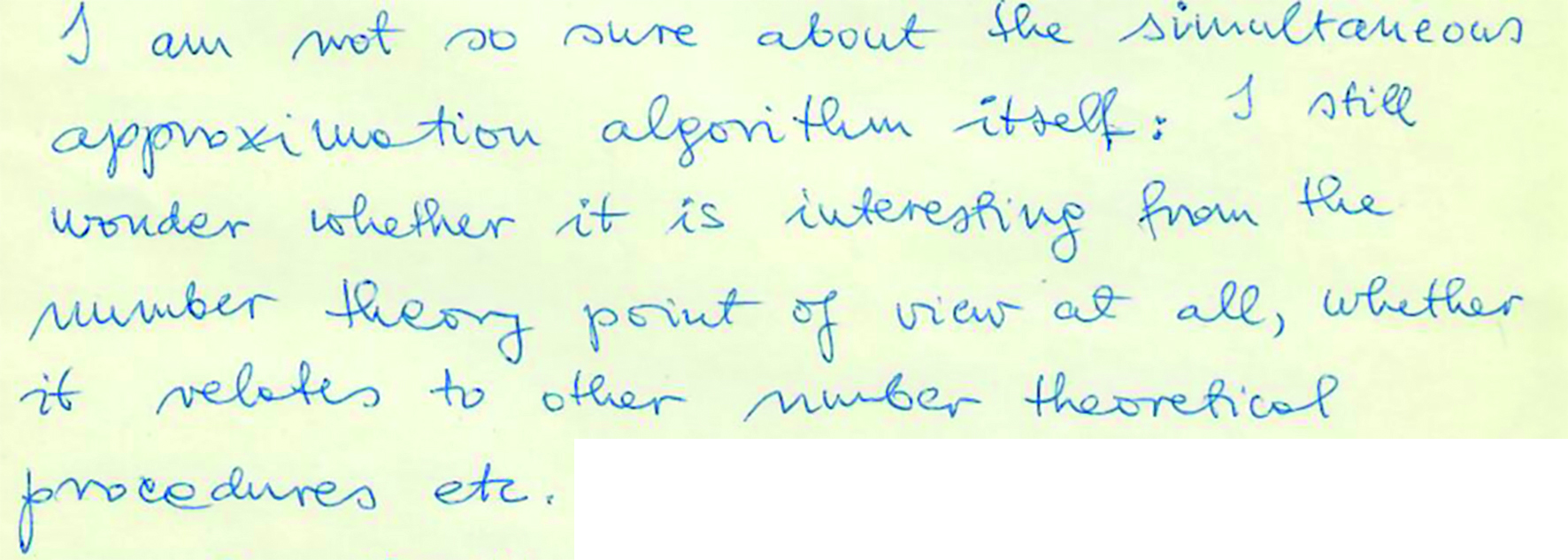}

\caption{Cutout from a Lovász letter}
\label{fig:GN11}
\end{figure}

\noindent
The basis reduction algorithm by L.~Lovász to solve a problem, that initially
looked like a technicality, had a significant impact on the book~\cite{GLS1988}
as outlined in Section~\ref{sec:GN11}. Its deep impact on other fields came
really unexpected, even for Lovász himself as can be inferred from his letter,
see Fig.~\ref{fig:GN11}.

We consider this as one of the occasional miracles in mathematics where
a result that was prompted by the desire to find an elegant solution for
a technical detail has consequences that are simply beyond imagination.

Lovász informed not only the coauthors Grötschel and Schrijver of his
book~\cite{GLS1988} about his achievement, but also Hendrik Lenstra.  Employing
tools from the geometry of numbers, Hendrik had (briefly before) made the
substantial discovery that integer programs (IPs) can be solved in polynomial
time when the dimension is fixed. Concerning this, he was in discussion with
Lovász who pointed out that some of the steps of Hendrik's IP-algorithm could
be improved, see~\cite{Lenstra1983}.

Hendrik got excited about the news because his brother Arjen was (together with
two fellow students) about to implement a method to factor univariate
polynomials over algebraic number fields. Zassenhaus had suggested to use the
Berlekamp--Hensel approach for this which, however, could be ``\emph{very, very
much exponential}'' according to Arjen. A few days after Lovász's letter had
arrived, Hendrik became convinced that the basis reduction algorithm implies
that there is a polynomial time algorithm for factorization in the ring
$\mathbb{Q}(X)$ of univariate polynomials over the rational numbers. At that
time this looked inconceivable as one did not (and still does not) know a
polynomial time algorithm for finding the factors of an integer. After working
out the details, Hendrik's observation turned out to be true.  The two Lenstra
brothers and Lovász combined their contributions and wrote the joint
paper~\cite{LenstraLenstraLovasz1982}. Believing that polynomial time factoring
of polynomials over the rational numbers (an unexpected result) is the most
important contribution of their work, they agreed to mention only this aspect
in the paper title. The full story of this cooperation is nicely described in
the article of I.~Smeets~\cite{Smeets2010}.

It turned out that basis reduction has applications that reach much further
than linear programming or polynomial factorization. It is beyond the scope of
this article to highlight here the wide range of applications of the basis
reduction algorithm which -- in contrast to the ellipsoid method -- is usable
in practice. We mention two concrete examples.

Odlyzko and te Riele~\cite{OdlyzkoteRiele1985} used the basis reduction
algorithm to disprove the Mertens conjecture, a conjecture standing in number
theory since~1897, which -- if true -- would have implied the Riemann
hypothesis. This disproof was surprising as there was extensive computational
evidence that the Mertens conjecture is true.

Lagarias and Odlyzko~\cite{LagariasOdlyzko1985} employed the lattice basis
reduction algorithm to launch a polynomial time attack on knapsack-based
public-key cryptosystems which made these cryptosystems unsafe.

The LLL algorithm, in fact, created a revolution in cryptography. It is known
that the widely used public-key schemes such as the RSA or elliptic-curve
cryptosystems can be defeated if Shor's quantum polynomial time factoring
algorithm can be implemented on a quantum computer. Many cryptographers are
convinced that certain lattice problems cannot be solved efficiently. Based on
this, some lattice-based constructions appear to be resistant to attack by both
classical and quantum computers. For surveys see Regev~\cite{Regev2006} or
Micciancio and Goldwasser~\cite{MicciancioGoldwasser2002}. The National
Institute of Standards and Technology (NIST) and other institutions are
currently preparing cryptography standards for the post-quantum era. The first
Quantum-Resistant Cryptographic Algorithms were announced by NIST in July~2022.
Lattices play a major role here, and lattice basis reduction algorithms have
become standard tools to test the security of cryptosystems.

Instead of attempting to comprehensively document the impact of Lovász's work
on basis reduction, we point to the book by Nguyen and
Vallée~\cite{NguyenVallee2010} entitled \emph{The LLL~Algorithm: Survey and
Applications} which consists of a collection of broad overviews of fields where
the LLL algorithm is employed. Chapters, written by specialists in the
respective fields, cover, for instance, applications in number theory,
Diophantine approximation, integer programming, cryptography, geometry of
provable security, inapproximability, and improvements of the LLL algorithm. A
reviewer of this book wrote: 

\begin{quote}
The LLL algorithm embodies the power of lattice reduction on a wide range of
problems in pure and applied fields [$\ldots$]  [and] the success of LLL attests to
the triumph of theory in computer science.
\end{quote}
Finally, the algorithm Lovász designed to find a reduced lattice basis is
usually called LLL algorithm, because it appeared in a paper written by three
authors whose last names starts with L. Of course, the Lenstra brothers do not
claim that it is their invention, they also attribute it to L.~Lovász. But LLL
algorithm has become the usually employed name of the algorithm.

\section{Cutting Planes and the Solution of Practical Applications}
\label{sec:GN13}

\begin{discussedreferences}
M. Grötschel, L. Lovász, A. Schrijver. \emph{Geometric
Algorithms and Combinatorial Optimization}. Springer, Berlin, 1988.
\end{discussedreferences}

\noindent 
László Lovász has, in addition to inventing beautiful theory, designed many
algorithms, concentrating particularly on polynomial time algorithms. The
theory and the algorithms Lovász developed had significant impact on
computational practice. Chapter 8 of \cite{GLS1988} ``Combinatorial
Optimization: A Tour d’Horizon'' is a highly condensed overview of the
applicational potential that arises from combinations of the many insights
provided by the ellipsoid method, the LLL algorithm, and further ideas. These
have contributed to the astonishing computational success stories that evolved
in the last thirty to forty years in combinatorial optimization. We sketch some
of these aspects.

In combinatorial optimization, a typical approach is, as indicated before, to
attack a problem by transforming it into a linear programming problem with
integer variables.

Take the traveling salesman problem, for instance. Given a complete graph~$G =
(V, E)$ on $n$~vertices and a distance~$c_{e}$ for every edge~$e \in E$, we
look for a Hamiltonian cycle (briefly: tour) of minimum length. If $H$ is a
tour, let $x^{H} \in \mathbb{R}^{E}$ be its \emph{incidence vector}, i.e., the
$e$-th component $x^{H}_{e}$ of $x^{H}$ is equal to $1$ if $e \in H$, otherwise
it is $0$. The traveling salesman polytope $\TSP(G)$ of $G$ is the convex hull
of all incidence vectors of tours in $G$. $\TSP(G)$ is a polytope in
$\mathbb{R}^{n(n-1)/2}$. To apply the linear programming approach, we now have
to find a linear inequality system, so that the integral solutions of the
linear program are exactly the incidence vectors of tours. Such linear programs
are called \emph{LP-relaxations}. Let $\delta(W)$ denote the set of edges
in~$E$ with one endvertex of~$e$ in~$W$ and the other in~$V\setminus W$, and
let $x(\delta(W))$ denote the sum over all variables~$x_{e}$ with $e \in
\delta(W)$. It is well known that the following linear program:
\begin{align*}
  0 \leq x_{e} \leq 1            & \quad \text{for all } e \in E\\
  x\bigl(\delta(\{w\})\bigr) = 2 & \quad \text{for all } w \in V\\
  x\bigl(\delta(W)\bigr) \geq 2  & \quad \text{for all } W \subseteq V \text{ with } 2 \leq |W| \leq |V| -2
\end{align*}
is an LP-relaxation of the $\TSP$. The third type of inequalities is called
\emph{subtour elimination constraints}.

Let us call the polytope defined by the linear system above $\TSPLP(G)$.  All
vertices of the traveling salesman polytope $\TSP(G)$ are vertices of
$\TSPLP(G)$. But $\TSPLP(G)$ has many nonintegral vertices as well. About
$2^{n}$ inequalities define $\TSPLP(G)$. This renders the straightforward
LP-solution approach hopeless. The facet complexity $\phi$ of $\TSPLP(G)$,
however, is small since the entries of every inequality or equation are only 0
or 1 and the right-hand sides are 0, 1, or 2. Thus the facet complexity of
$\TSPLP(G)$ is linear in the number of variables $|E| = n(n-1)/2$. Due to the
oracle-polynomial time equivalence of strong separation and strong
optimization, linear programs over $\TSGLP(G)$ can be solved in polynomial time
-- provided, given a vector $y \in \mathbb{Q}^{E}$, one can find a fast
separation algorithm for the subtour elimination constraints.

This can in fact be done, as was observed by Hong~\cite{Hong1972}. One
assigns the value $y_{e}$ to every edge $e \in E$ as a capacity
and computes (this can be done quickly) a minimum nonempty cut $\delta(W^*)$ in
this capacitated graph $G = (V, E)$. If $y(\delta(W^*)) < 2$, a violated
inequality is found, otherwise y satisfies all subtour elimination
constraints. This is an example of a linear program appearing in many
practical applications with an exponential number of inequalities that,
nevertheless, can be solved in polynomial time. An optimal solution of a
linear program over $\TSPLP(G)$ is usually nonintegral but provides a very
good lower bound on the optimum $\TSP$-value in practice. Finding a
provably optimal solution needs additional effort, though.

In 1954 Dantzig, Fulkerson, and Johnson~\cite{DantzigFulkersonJohnson1954}
proposed in a seminal paper to solve combinatorial optimization problems such
as the traveling salesman problem by starting with some LP-relaxation, checking
whether the optimum solution~$y$ is the incidence vector of a tour (in this
case the problem is solved), and if not searching for inequalities valid for
TSP(G) that are violated by~$y$, adding these to the current LP as
\emph{cutting planes}, and to continue. This was one of the first proposals to
solve linear and integer programs using cutting planes in an iterative process.
The cutting plane search in this case was done manually, the LPs were solved by
the simplex method. Four years later Gomory~\cite{Gomory1958} invented an
automatic cutting plane generation scheme (called \emph{Gomory cuts}) for which
he could prove finite termination. This looked like a promising approach to
solve integer programs.

However, the computer implementations of this and related approaches in the
1960s and 1970s were not successful in practice. Moreover, theoretical results
of Chvátal~\cite{Chvatal1973} and others showed that there are series of
examples for which the number of cutting plane additions cannot be effectively
bounded. Hoping for unimportance of these negative aspects in real-world
applications, the idea came up in the 1970s to study combinatorial optimization
problems of practical relevance and to look for cutting planes that define
facets of the investigated polytopes. These are cuts that cut as deep as
possible. The first implementations employing a combination of manual and
heuristic searches for facet defining cutting planes at the end of the~1970s
indicated practical success. Soon after, the ellipsoid method theory with the
principle of polynomial time equivalence of optimization and separation was
developed and demonstrated that this approach is a viable idea, and that linear
optimization over exponentially large systems of linear inequalities is
possible in polynomial time -- at least theoretically.

Despite serious attempts, no implementation of the ellipsoid method has shown
satisfactory numerical performance in computational practice. By replacing it
with new implementations of the dual simplex algorithm, the theoretical
polynomial time termination is lost, but astonishing computational results were
achieved by many researchers in combinatorial optimization. Of course, lots of
additional features (such as presolve techniques, heuristic primal and dual
searches, branch and bound, robust numerics, etc.) were implemented as well.
The new insights gave a significant push to the theoretical and applied part of
combinatorial optimization. Problems with many industrial applications such as
linear ordering; set partitioning and packing; knapsack; clustering; various
types of matching; connectivity; path, flow and other network problems; max
cut; unconstrained Boolean quadratic programming; stable sets; several
variations of coloring; and vehicle and passenger routing could be solved for
instances of practically relevant sizes. The discovery of new classes of facets
and fast separation procedures (exact and heuristic) has been an important
ingredient of this solution methodology. To indicate at least one example of
practically useful separation algorithms we mention the
paper~\cite{PadbergRao1982} of Padberg and Rao that describes sophisticated and
fast separation algorithms for various ramifications of the matching polytope.
A large number of separation algorithms are, of course, described in the
book~\cite{GLS1988}.

This research activity goes on and brings application relevant instances of
many \cNP-hard combinatorial optimization problems to the realm of
practical solvability. For the traveling salesman problem, for example, the
``solvability world record'' was 42~cities in 1954, it went to~120 in~1977,
2392 in~1987, and in~2017 a TSP with 109,399~cities could be solved to
optimality, see the Webpage of Bill Cook~\cite{Cook2023}, his
book~\cite{Cook2012}, and the book~\cite{ApplegateBixbyChvatalCook2007} by
Applegate, Bixby, Chvátal, and Cook for comprehensive information. The solution
process includes linear programming technology (its theory and implementation)
that is able to prove, for example, that a vector in dimension $10^{10}$
satisfies more than $2^{100{,}000}$ constraints and is optimal for this system.
This is really breathtaking.

\begin{figure}[t]
\centering
\includegraphics[width=9cm]{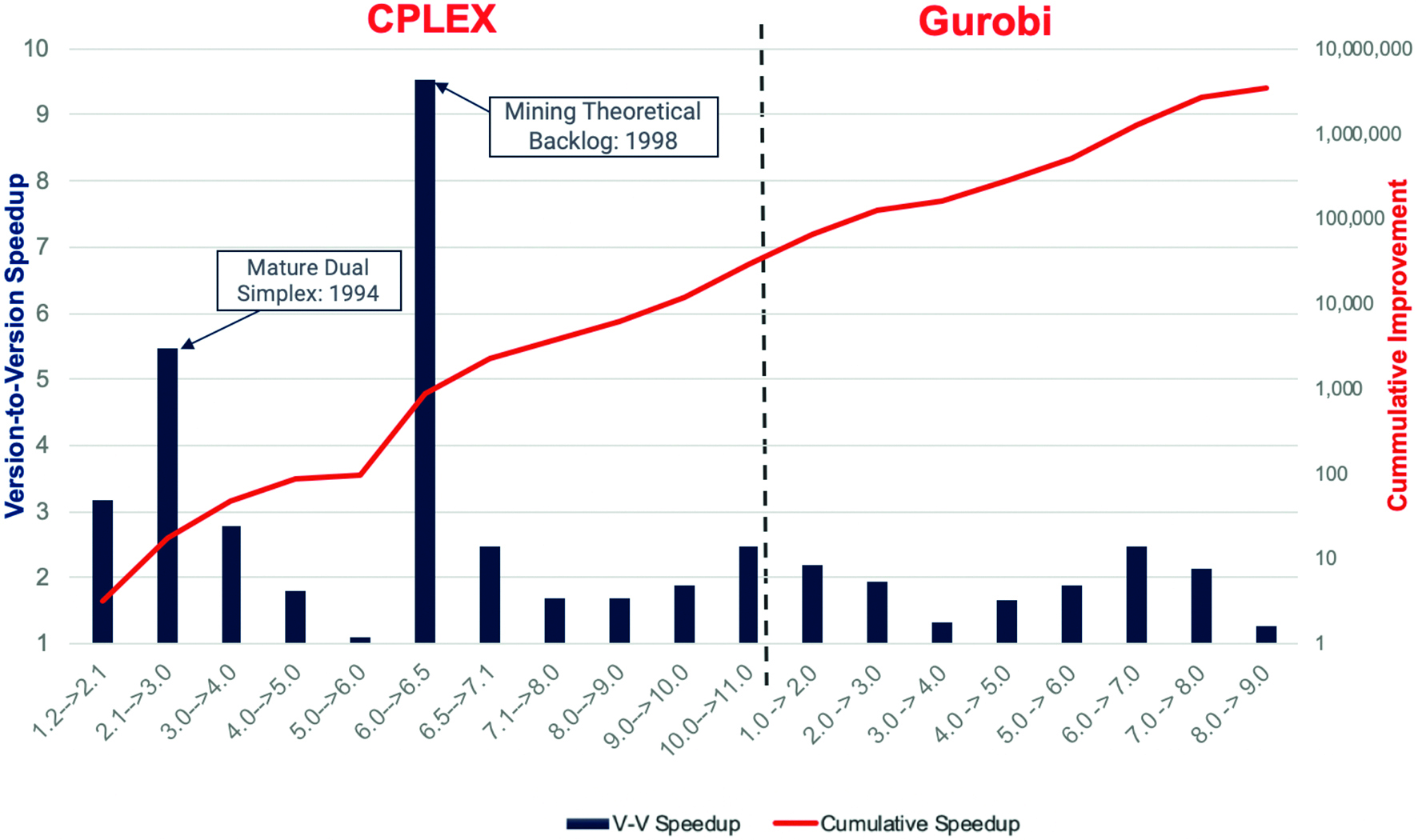}

\caption{MIP-code performance 1990--2019 (courtesy Robert E. Bixby)}
\label{fig:GN12}
\end{figure}

The success stories indicated above, and the theoretical and practical lessons
learned from these began to be harvested and improved by the developers of
commercial optimization software in the 1990s. One reason for this is that many
mixed-integer optimization problems (MIPs) occurring in industry contain
subproblems that are combinatorial optimization problems for which large
classes of facet-defining inequalities have been discovered. Efficient
separation algorithms for these inequalities were successfully added to the
existing MIP-codes.  The graphic in Fig.~\ref{fig:GN12}, presented with the
permission of Bob Bixby, shows the development of the commercial mixed integer
programming codes CPLEX and Gurobi in the 30~years from~1990 to~2019. The large
bar (pointed at by ``Mining Theoretical Backlog'') shows an almost tenfold
speedup that is obtained from one version of the code to the next in which
cutting plane technology (including a fresh implementation of Gomory cuts) was
introduced together with various supporting features.  The overall message is
that the MIP technology in~2019 runs 3.5~million times faster than the codes
of~1990. That speedup is due to mathematical and implementation improvements
and is independent of the hardware speedup during this period. This is real
progress indeed. Cutting plane technology contributed to it significantly.

\section{Computing Optimal Stable Sets and Colorings in Perfect Graphs}
\label{sec:GN14}

\begin{discussedreferences}
M. Grötschel, L. Lovász, A. Schrijver. Polynomial algorithms for perfect
graphs. \emph{Annals of Discrete Math}. 21:325--256, 1984.

M. Grötschel, L. Lovász, A. Schrijver. Relaxations of vertex packing.
\emph{J. Combin. Theory} B 40:330--343, 1986.
\end{discussedreferences}

\noindent
The extension of the ellipsoid method to convex bodies outlined in
Section~\ref{sec:GN10} was driven by the hope that one could solve the stable
set and the coloring problem in perfect graphs in polynomial time with this
methodology, see Section~\ref{sec:GN7}. The successful attempt is presented in
the articles~\cite{GLS1981, GLS1984}, and~\cite{GLS1986}. We describe the
stable set case.

For a graph $G=(V, E)$ and a stable set $S \subseteq  V$, one can define the
incidence vector~$x^{S}$ in $\mathbb{R}^{V}$ as
follows: the $i$-th component $x^{S}_{i}$ of
$x^{S}$ is equal to 1 if the vertex $i \in V$ is an element of
$S$, and it is 0 otherwise. The stable set polytope of $G$ is the convex
hull of all incidence vectors of stable sets $S$ of $G$, i.e.,
\[
\STAB(G) \coloneqq \conv\bigl\{x^{S} \in \mathbb{R}^{V} \mid
S\subseteq V \text{ stable set}\bigr\}.
\]
Let $w: V\rightarrow \mathbb{Q}$ be any weighting of the vertices of $G$ (we
may assume that all weights are positive) and denote the largest weight of a
stable set in $G$ by $\alpha(G,w)$. Then $\alpha(G,w)$ is the maximum value of
the linear function $w^{T}x$ for $x \in \STAB(G)$, in other words,
$\alpha(G,w)$ can be computed by solving a linear program over $\STAB(G)$. For
this observation to be of any use, we have to find inequalities defining
$\STAB(G)$.  Consider the polytope defined by
\[
\begin{split}
  \QSTAB(G) \coloneqq \bigl\{x \in \mathbb{R}^{V}   \mid & x_{i} \geq 0
\quad\forall i \in V, \  
x_{i} + x_{j} \leq 1 \quad\forall ij \in E, 
\\
& x(Q) \leq 1 \quad \forall Q \subseteq  V \text{ clique}\bigr\},
\end{split}
\]
where $x(Q)$ denotes the sum of all $x_{i}, i \in Q$. The corresponding
inequality is called \emph{clique constraint}. Since the intersection of a
clique and a stable set contains at most one vertex, all clique constraints are
satisfied by all incidence vectors of stable sets. This implies $\STAB(G)
\subseteq  \QSTAB(G)$ and optimizing over $\QSTAB(G)$ is an LP-relaxation of
the stable set problem.

The stable set problem is \cNP-hard. Therefore, solving linear programs
over $\STAB(G)$ is \cNP-hard as well. For some combinatorial
optimization problems, their natural LP-relaxation is solvable in polynomial
time. A sobering observation is that, for general graphs, solving linear
programs over $\QSTAB(G)$ is also \cNP-hard.  So, in general, nothing is
gained algorithmically. For perfect graphs, though, this approach combined with
a tighter relaxation delivers the desired result.

Lovász's Shannon capacity article~\cite{Lovasz1979a} suggests studying a
different relaxation of the stable set problem.

Let $(u_i \mid i \in V), u_i \in \mathbb{R}^N$, be any orthonormal
representation of $G$ and let $c \in \mathbb{R}^N$ with $\|c\| = 1$. Then for
any stable set $S \subseteq V$, the vectors $u_i, i \in S$, are mutually
orthogonal and hence,
\[
  \sum_{i \in S} (c^T u_i)^2 \leq 1.
\]
Since $\sum_{i \in V} (c^T u_i)^2 x_i^S = \sum_{i \in S} (c^T u_i)^2$, we see
that the inequality
\begin{equation}
\tag{ORC}
\sum_{i \in V} (c^T u_i)^2 x_i \leq 1
\end{equation}
holds for the incidence vector $x^S \in \mathbb{R}^V$ of any stable $S$ set of
nodes of $G$. Thus, (ORC) is a valid inequality for $\STAB(G)$ for any
orthonormal representation $(u_i \mid i \in V)$ of~$G$, where $u_i \in
\mathbb{R}^N$, and any unit vector $c \in \mathbb{R}^N$. We shall call (ORC)
the \emph{orthonormal representation constraints} for $\STAB(G)$.

Utilizing these inequalities, the following set was introduced
in~\cite{GLS1986}. For any graph $G=(V,E)$ let
\begin{multline*}
\TH(G) \coloneqq \bigl\{x \in \mathbb{R}^{V} \mid x_{i} \geq 0 \quad \forall
i \in V, \\
\text{ and $x$ satisfies all orthonormal representation constraints}\bigr\}.
\end{multline*}
$\TH(G)$ is the solution set of infinitely many linear inequalities and thus a
convex set. Since for every clique $Q$, its clique constraint appears as an
orthonormal representation constraint (given a clique $Q\subseteq V$, let
$\{u_{i} \mid i \in V\setminus Q\} \cup \{c\}$ be mutually orthogonal unit
vectors and set $u_{j} = c$ for$ j \in Q$) and every incidence vector of a
stable set satisfies all such inequalities, we obtain:
\[
\STAB(G) \subseteq  \TH(G) \subseteq  \QSTAB(G).
\]
An important fact is, that the Lovász theta function $\vartheta(G,w)$ introduced
in Section~\ref{sec:GN8} can also be characterized as follows:
\[
\vartheta(G,w) = \max \bigl\{w^{T}x \mid x \in \TH(G)\bigr\}.
\]
$\TH(G)$ is contained in the unit ball, and it is easy to find the center of a
ball contained in the interior of $\TH(G)$. Thus, $\TH(G)$ is a convex body
satisfying the assumptions required for the oracle-polynomial time equivalence
of weak optimization, separation, and membership. The desired result is, of
course, the following:

\begin{nntheorem}
The weak optimization problem for $\TH(G)$ is solvable in polynomial time for
any graph $G = (V,E)$.
\end{nntheorem}

Lovász, see \cite{Lovasz1979a}, established several characterizations for his
$\vartheta$-function.  They can be used in various ways to prove this theorem.
One proof, worked out in detail in~\cite{GLS1981} and \cite{GLS1988}, is based
on the following characterization:
\begin{align*}
& \vartheta(G,w) =  \max \{\bar{w}^{T}B\bar{w} \mid B \in \mathcal{K} \}, \\
& \text{where }
\mathcal{K} \coloneqq \{B \in \mathbb{R}^{V \times V} \mid B \in \mathcal{D} \cap \mathcal{M} \text{ and } \tr(B)=1\}.
\end{align*}
Above, $\mathcal{D}$ is the set of positive semidefinite matrices,
$\mathcal{M}$ the set of symmetric matrices $B$ that satisfy $b_{ij}=0$
whenever $ij$ is an edge in $G$, and $\bar{w}$ denotes the vector whose entries
are the square roots of the values $w_{i}, i \in V$. The main part of the proof
consists in showing that the weak membership problem for $\mathcal{K}$ can be solved in
polynomial time, and the core of this proof is established by showing whether a
symmetric matrix is positive definite.

A by-product of the proof is the first polynomial time algorithm for
optimization problems containing positive semidefinite constraints, a major
result that led to considerable follow-up research such as the design of
polynomial time interior point (and other) algorithms for semidefinite
programming.

Another way to establish the above theorem is by utilizing the following fact:
\[
  \vartheta(G,w) = \min \bigl\{ \Lambda (A+W) \mid A \in \mathcal{M}^{\perp}\bigr\},
\]
where $\Lambda$ denotes the largest eigenvalue, $\mathcal{M}^{\perp}$ the
orthogonal complement of $\mathcal{M}$, and $W$ the symmetric $V\times
V$-matrix whose entries are the square roots of $w_{i}w_{j}$.  $\Lambda(A+W)$
is a convex function that ranges over a linear space, and thus, we can obtain
$\vartheta(G,w)$ via an unconstrained convex function optimization problem in
polynomial time.

A third way to prove the theorem was demonstrated
in~\cite{LovaszSchrijver1991}, and this approach turned out to be one of the
starting points for a generalization of this technique. Lovász and Schrijver
developed in this article a general lift-and-project method that constructs
higher-dimensional polyhedra (or, in some cases, convex sets) whose projection
approximates the convex hull of 0-1 valued solutions of a system of linear
inequalities. An important feature of these approximations is that one can
optimize any linear objective function over them in polynomial time.
Lift-and-project methods have been extended in many directions and are still an
area of intensive research. The recent (not even exhaustive) survey by Fawzi,
Gouveia, Parrilo, Saunderson, and
Thomas~\cite{FawziGouveiaParriloSaundersonThomas2022} discusses the
contributions of almost one hundred articles and illustrates the richness of
this topic by presenting examples from many different areas of mathematics and
its applications.

We refrain from describing the technically challenging details of this
lift-and-project technique and return to stable sets in perfect graphs.

A combination results of Fulkerson~\cite{Fulkerson1972} and
Chvátal~\cite{Chvatal1975} yields:

\begin{nntheorem}
$\STAB(G) = \QSTAG(G)$ if and only if $G$ is perfect.
\end{nntheorem}

And since we already know that $\STAB(G) \subseteq  \TH(G) \subseteq
\QSTAB(G)$ holds, we obtain:

\begin{nncorollary}
$\STAB(G) = \TH(G) = \QSTAG(G)$ if and only if G is perfect.
\end{nncorollary}

Since the weak optimization problem for $\TH(G)$ can be solved in polynomial
time and since, in case $G$ is perfect, $\TH(G)$ is a well-described
polyhedron, the strong optimization problem for $\TH(G)$ can be solved in
polynomial time. This yields the desired result:

\begin{nntheorem}
The stable set problem can be solved in polynomial time for perfect graphs.
\end{nntheorem}

We can now employ the fact that, if a linear program can be solved in
polynomial time, the dual linear program can also be solved in polynomial time,
see Section~\ref{sec:GN11}. By proving that, in this case, an optimum basic
solution of the dual program can be transformed in polynomial time into an
integral optimum basic solution one can find an optimum solution of the
weighted clique covering problem. Since the cliques of a graph~$G$ are the
stable sets of the complementary graph~$\bar{G}$ of $G$ and the colorings
of~$G$ are the clique covering of $\bar{G}$, we can conclude:

\begin{nntheorem}
For perfect graphs, the stable set, the clique, the coloring, and the clique
covering problem can be solved in polynomial time. This also holds for the
weighted versions of these problems.
\end{nntheorem}

\section{Submodular Functions}
\label{sec:GN15}

\begin{discussedreferences}
L. Lovász. Submodular functions and convexity. In \emph{Mathematical
Programming: The State of the Art} (eds. A. Bachem, M. Grötschel, B.
Korte), Springer, pages 235--257, 1983.
\end{discussedreferences}

\noindent
Let $E$ be a finite set. A function $f: 2^{E} \rightarrow \mathbb{R}$ is called
\emph{submodular} on $2^{E}$ (the power set of $E$) if
\[
f(S\cap T) + f(S\cup T) \leq f(S) + f(T) \text{ for all } S, T \subseteq  E.
\]
Submodular functions play an important role in lattice theory, geometry, graph
theory, and particularly, in matroid theory and matroidal optimization
problems. The rank function of a matroid, for example, is submodular as well as
the capacity function of the cuts in directed and undirected graphs.

Two polyhedra can be associated with a submodular function $f:
2^{E} \rightarrow \mathbb{R}$ in a natural way
\begin{align*}
& P_{f} \coloneqq \bigl\{x \in \mathbb{R}^{E} \mid x(F) \leq f(F) \text{ for all } F \subseteq  E, \ x \geq 0\bigr\},
\\
& EP_{f} \coloneqq \bigl\{x \in \mathbb{R}^{E} \mid x(F) \leq f(F)
\text{ for all } F \subseteq  E\bigr\}.
\end{align*}
$P_{f}$ is called the \emph{polymatroid} associated with the submodular
function~$f$, $EP_{f}$ the \emph{extended polymatroid} associated with~$f$. A deep
theorem of Edmonds~\cite{Edmonds1970} states that if $f$ and $g$ are two
integer valued submodular functions then all vertices of $P_{f}\cap P_{g}$ as well
as all vertices of $EP_{f}\cap EP_{g}$ are integral. This theorem contains a large
number of integrality results in polyhedral combinatorics; it particularly
generalizes the matroid intersection theorem.

To address algorithmic questions concerning the structures introduced above, we
assume that a submodular function~$f$ is given by an oracle that returns the
value~$f(S)$ for every query $S \subseteq  E$. We also assume that we know an
upper bound $\beta$ on the encoding length of the output of the oracle. With
these assumptions we define the encoding length of the submodular function as
$|E| + \beta$.

It is well-known that, for any nonnegative linear objective function, the
greedy algorithm finds an optimum vertex of $EP_{f}$ in oracle-polynomial time,
and that this vertex is integral provided the submodular function~$f$ is
integer valued. Optimizing over polymatroids or the intersections of two
polymatroids or the intersections of two extended polymatroids and finding
integral optima is more complicated and needs careful analysis. The most
important algorithmic problem in this context is:

\begin{SFM}
Given a submodular function $f: 2^{E} \rightarrow \mathbb{Q}$, find a set~$S
\subseteq  E$ minimizing $f$.
\end{SFM}

Lovász has built in~\cite{Lovasz1983} a bridge between submodularity and
convexity by showing that submodular functions are discrete analogues of convex
functions and has thus provided the key to the algorithmic solution of the
submodular function minimization problem. The link is established as follows.

Let $f: 2^{E} \rightarrow \mathbb{R}$ be any set function. For every subset $T
\subseteq  E$, let $x^{T}$ be its incidence vector and set
\[
  \hat{f} (x^{T}) \coloneqq f(T).
\]
This way $\hat{f}$ is defined on all 0/1-vectors. Note that every nonzero
nonnegative vector~$y \in \mathbb{R}^{E}$ can be expressed uniquely as
\begin{multline*}
y = \lambda_{1}x^{T_1} +
\lambda_{2}x^{T_2} + \ldots  + \lambda_{k}x^{T_k}, 
\\
\text{such that}\quad 
\lambda_{i} > 0, i = 1,\ldots,k 
\quad\text{and}\quad
\emptyset \neq T_{1} \subset T_{2} \subset \ldots  \subset T_{k} \subseteq  E.
\end{multline*}
Then
\[
\hat{f}(y) \coloneqq \lambda_{1}f(T_{1})+
\lambda_{2}f(T_{2}) + \ldots +
\lambda_{k}f(T_{k})
\]
is a well-defined extension of the set function $f$ (called \emph{Lovász
extension of $f$}) to the nonnegative orthant. Lovász proved
in~\cite{Lovasz1983}:

\begin{nntheorem}
Let $f: 2^{E} \rightarrow \mathbb{R}$ be any set function and $\hat{f}$ its
extension to nonnegative vectors. Then $\hat{f}$ is convex if and only if $f$
is submodular.
\end{nntheorem}

\begin{nnlemma}
Let $f: 2^{E} \rightarrow \mathbb{R}$ be set function with $f(\emptyset) =
0$. Then
\[
\min \bigl\{f(S) \mid S \subseteq  E\bigr\} = \min \bigl\{\hat{f}(x) \mid x \in [0, 1]^{E}\bigr\}.
\]
\end{nnlemma}

Thus, instead of minimizing a set function $f$ over $E$, it suffices to
minimize its Lovász extension $\hat{f}$ over the unit hypercube. We observe
that $\hat{f}(x)$ can be evaluated in oracle-polynomial time using the oracle
defining $f$ and that, if $f$ is submodular, then $\hat{f}$ is convex. We know
already from Section~\ref{sec:GN10} that convex functions can be minimized in
oracle-polynomial time. (The assumption $f(\emptyset) = 0$ is irrelevant, if
necessary, we can replace $f$ by the function $f-f(\emptyset)$.) This yields:

\begin{nntheorem}
Let $f: 2^{E} \rightarrow \mathbb{Q}$ be a submodular function. Then a subset
$S$ of $E$ minimizing $f$ can be found in oracle polynomial time.
\end{nntheorem}

This theorem implies the polynomial time solvability of many combinatorial
optimization problems, including the computation of a minimum capacity cut in a
graph. It has various ramifications such as solvability in strongly polynomial
time, as outlined in~\cite{Lovasz1983} and~\cite{GLS1988}.

The running time of the polynomial time algorithm sketched above makes it,
however, infeasible for practical use. New and better polynomial time
algorithms, not employing the ellipsoid method, have been devised by
Schrijver~\cite{Schrijver2000} and Iwata, Fleischer, and
Fujishige~\cite{IwataFleischerFujishige2001}.

\section{Volume Computation}
\label{sec:GN16}

\begin{discussedreferences}
L. Lovász. How to compute the volume? \emph{Jber. d. Dt.
Math.-Vereinigung}, Jubi\-läums\-tagung 1990, B.\,G. Teubner, Stuttgart,
pages 138--151, 1992.
\end{discussedreferences}

\noindent
Since the convergence of all versions of the ellipsoid method depends on
sequentially shrinking the volume of an ellipsoid containing the given convex
body $K$, it is tempting to ask whether the algorithm can be tuned to provide a
reasonable estimate of the volume of~$K$. The key idea in this context is, of
course, to come up with an algorithmic version of the Löwner--John theorem,
that states, that, for a convex body~$K$ in $\mathbb{R}^{n}$, there exists a
unique ellipsoid~$E$ of minimal volume containing $K$; moreover, $K$ contains the
ellipsoid obtained from $E$ by shrinking it from its center by a factor of $n$.
In formulas, let $E(A, a) \coloneqq \{ x \in \mathbb{R}^{n} \mid
(x-a)^{T}A^{-1}(x-a) \leq 1\}$ denote the ellipsoid defined by a positive
definite matrix $A$ with center $a \in \mathbb{R}^{n}$ then the Löwner--John
theorem states
\[
E(n^{-2}A, a) \subseteq  K \subseteq  E(A, a),
\]
if $E(A, a)$ is the Löwner--John ellipsoid $E$ of $K$. Algorithmically, the
following could be achieved in the Grötschel--Lovász--Schrijver
book~\cite{GLS1988}.

\begin{nntheorem}
There exists an oracle-polynomial time algorithm that finds, for
any convex body $K$ given by the space dimension $n$, a weak separation
oracle and two real numbers $r$ and $R$ with the property that $K$ is
contained in the ball of radius $R$ around the origin and contains a ball
of radius $r$, an ellipsoid $E(A, a)$ such that
\[
  E\Bigl(\frac{1}{n(n+1)^{2}} A, a\Bigr) \subseteq  K \subseteq  E(A, a).
\]
\end{nntheorem}

With more effort and making additional assumptions such as central symmetry or
requiring that a system of defining linear inequalities is explicitly given (in
the polytopal case), the factor $1/(n(n+1)^{2}$ in front of the matrix~$A$ above
can be slightly improved, but not fundamentally. If one declares the volume of
the interior ellipsoid as an approximation of the volume of $K$, the relative
error turns out to be $2^{n}n^{3n/2}$, which appears to be outrageously bad.

Surprisingly, the error is not as bad as it looks since subsequently
Elekes~\cite{Elekes1986} and others proved that no oracle-polynomial time
algorithm can compute, for a convex body $K$ as given above, the volume of $K$
with a much better relative error. We quote a result of Bárány and
Füredi~\cite{BaranyFuredi1986}.

\begin{nntheorem}
Consider a polynomial time algorithm which assigns to every convex body~$K$
given by a membership oracle an upper bound~$w(K)$ on its volume $\vol(K)$.
Then there is a constant $c>0$ such that in every dimension $n$ there exists a
convex body~$K$ for which $w(K) > n^{cn}\vol(K)$.
\end{nntheorem}

Following up, various authors proved more negative results on the deterministic
approximation of the volume, width, diameter and other convexity parameters.

These negative results fueled the investigation of stochastic approaches to
estimate the volume of a convex body. Instead of giving a deterministic
guarantee, one could try to calculate a number that is close to the true value
of the volume with high probability employing a randomized algorithm.

A side remark: Khachiyan~\cite{Khachiyan1993} and Lawrence~\cite{Lawrence1991}
proved that, for every dimension $n$, one can construct systems of rational
inequalities defining polytopes $P$ so that the encoding length of the rational
number $p/q$ representing the true volume of $P$ requires a number of digits
that is exponential in the encoding length of the inequality system. Hence,
exact volumes of convex bodies cannot be computed in polynomial time since
specifying the exact volume requires exponential space.

A fundamental breakthrough was achieved in Dyer, Frieze, and
Kannan~\cite{DyerFriezeKannan1991} who provided a randomized polynomial time
approximation scheme for the volume approximation problem where $K$ is given by
a membership oracle. The ingredients of their algorithm are a multiphase
Monte-Carlo algorithm (using the so-called product estimator) to reduce volume
computation to sampling, the utilization of Markov chain techniques for
sampling, and the use of the conductance bound on the mixing time, due to
Jerrum and Sinclair~\cite{JerrumSinclair1988}.  The running time of the
algorithm is roughly~$O(n^{23})$ which is truly prohibitive. The exponent~23 of
$n$ was subsequently reduced considerably by adding further techniques and
improved estimates to the toolbox of randomized algorithms, including rapid
mixing, harmonic functions, connection to the heat kernel, isoperimetric
inequalities, discrete forms of Cheeger inequality, and many more.

Lovász played an important role in the exponent shrinking race. For example,
the exponent went down to~16 (Lovász and
Simonovits~\cite{LovaszSimonovits1990}), to~10 (Lovász~\cite{Lovasz1992}), to~8
(Dyer and Frieze~\cite{DyerFrieze1992}), to~7 (Lovász and
Simonovits~\cite{LovaszSimonovits1993}), to~5 (Kannan, Lovász, and
Simonovits~\cite{KannanLovaszSimonovits1997}), and to~4 (Lovász and
Vempala~\cite{LovaszVempala2006}). A nice survey of the many tricky issues in
designing randomized algorithms for volume computation and their analysis is
the article by Simonovits~\cite{Simonovits2003}.

The race for better algorithms has not stopped. On September~3, 2022, the new
record was published on arXiv by Jia, Laddha, Lee, and
Vempala~\cite{JiaLaddhaAditiLeeVempala2022}. The authors show that the volume
of a convex body in $\mathbb{R}^{n}$ defined by a membership oracle can be
computed to within relative error $\epsilon$ using $\tilde{O}(n^{3}\psi^{2} +
n^{3}/\epsilon^{2})$ oracle queries, where $\psi$ is the KLS constant. With the
current bound of $\psi = \tilde{O}(1)$, this gives an
$\tilde{O}(n^{3}/\epsilon^{2})$ algorithm, improving on the Lovász--Vempala
$\tilde{O}(n^{4}/\epsilon^{2})$ algorithm.

\section{Analysis, Algebra, and Graph Limits}
\label{sec:GN17}

\begin{discussedreferences}
L. Lovász, \emph{Large Networks and Graph Limits}. American Mathematical
Society, 2012.
\end{discussedreferences}

\begin{figure}[b]
\centering

\setlength{\fboxsep}{0pt}
\setlength{\fboxrule}{.25pt}
\fbox{\includegraphics[width=2.5cm]{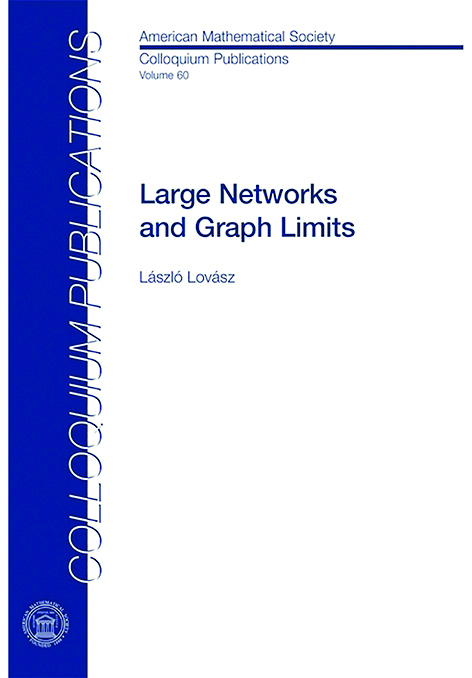}}

\caption{Lovász's Graph Limits book}
\label{fig:GN13}
\end{figure}

\noindent We have already indicated that many of the results mentioned in our
article seem to be of permanent importance and are used again and again: the
Lovász Local Lemma, algorithmic consequences of the ellipsoid method,
topological combinatorics, and the LLL algorithm, to name just few. Very
recently Lovász's mathematics culminated in a topic that somehow combines this
into an all-in-one subject: like a late symphony of a grand composer displaying
the experience of the master and an echo of his/her life. We believe that this
happened with the subject of graph limits founded and developed by Lovász with
co-authors and students in the last 15~years. Here is a brief sketch of this
fascinating development.

We have seen in Section~\ref{sec:GN2} that the homomorphism function
$\hom(F,G)$ and the Lovász vector $L(G)$ are determining every graph $G$ up to
an isomorphism. With a proper scaling this leads to the notion of
\emph{homomorphism density} $t(F,G)$, which is the probability that a random
mapping between sets of vertices of $F$ and $G$ is a homomorphism: \smash{$t(F,G) =
\frac{\hom{(F,G)}}{{v(G)}^{v(F)}}$} where $v(G)$ denotes the number of vertices
of graph $G$.

This definition is close to the sampling density and one motivation for
introducing it. One can observe that homomorphism densities do not determine a
graph up to an isomorphism but up to a ``blowing up of vertices''. (This is a
procedure by which vertices are replaced by a certain number of twin copies.)
It is perhaps more important that one can then define \emph{convergence of a
sequence} of finite graphs $G_{1},G_{2},\ldots, G_{n}, \ldots$ as the
convergence of homomorphism densities $t(F,G_{n})$ for every graph $F$.  This
convergence concept (and various other notions of convergence) were introduced
and investigated in the article \cite{BorgsChayesLovaszSosVesztergombi2008} of
C.~Borgs, J.\,T.~Chayes, L.~Lov\'asz, V.\,T.~S\'os, and K.~Vesztergombi.

Hence, a sequence of graphs converges if, for every $F$, all homomorphism
densities (or $F$-sampling densities) converge. Does this convergence have a
real (geometrical) meaning? Are there limit graphs or, perhaps, other limit
objects?

It appears that these questions have non-trivial yet positive answers and these
were the starting point of a very rich and interesting area.  In fact, they
generated a whole new theory. Here is a sample of some of the results.

L. Lovász and B. Szegedy proved the following in~\cite{LovaszSzegedy2006}:

\begin{nntheorem}
A sequence of graphs (with unbounded size) is converging if and
only if it converges to a symmetric measurable function
$W: [0,1]^{2} \rightarrow [0,1]$.
Moreover, up to a measurable bijection, such a function W is
uniquely determined.
\end{nntheorem}

\begin{figure}
\centering

\includegraphics[height=2.5cm]{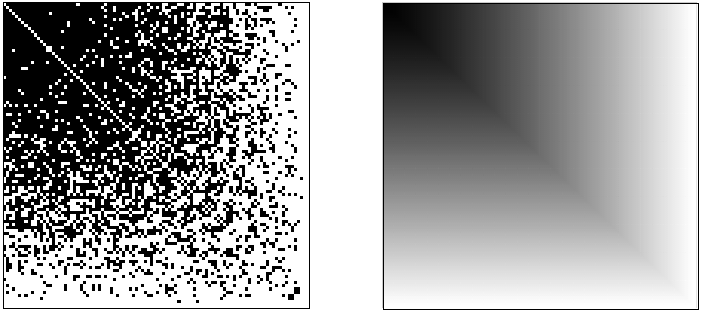}

\bigskip

\includegraphics[height=2.5cm]{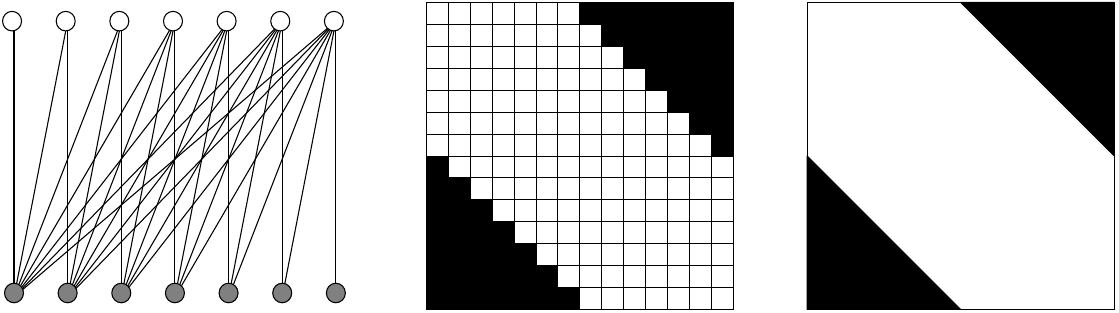}

\bigskip

\includegraphics[height=2.5cm]{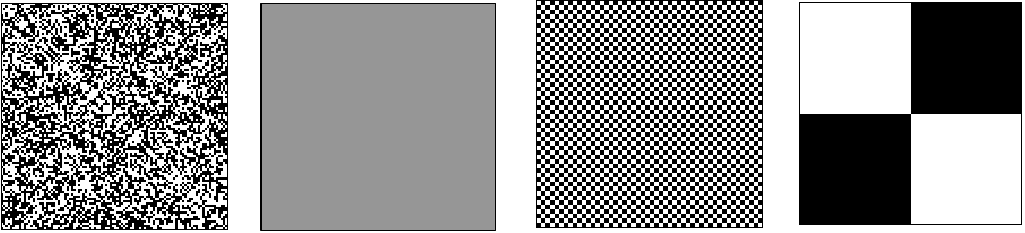}

\caption{Samples of graphons}
\label{fig:GN14}
\end{figure}

Explicitly, this means that for every graph
$F = (V,E)$ the homomorphism densities $t(F, G_{n})$ are
converging to:
\[
  t(F,W) = \int\limits_{[0,1]^{V}} \prod_{ij \in E} W(x_{i},x_{j}) \prod_{i \in V} dx_{i}
\]
Such functions $W$ are called \emph{graphons}. Graphon is a very intuitive
notion and the convergence of a graph sequence to a graphon looks like a
movie. It leads to ``pixel'' pictures like those on samples shown in
Figure~\ref{fig:GN14} (taken from Lovász`s book~\cite{Lovasz2012}).

The first row of Fig.~\ref{fig:GN14} shows on the left a randomly grown uniform
attachment graph with 100~nodes, and on the right a (continuous) function
approximating it. The picture on the right side is a grayscale image of the
function $U(x, y) = 1-\max(x,y)$.  The second row of Fig.~\ref{fig:GN14}
indicates the construction of the graphon for the ``halfgraph'' (the graph on
the left side). The bottom part indicates the influence of ordering and the
regularity Lemma in its simplest form. Note that the sequence of random graphs
is converging to a graphon $W$ that is a constant function. It is important
that the same is true for ``quasirandom graphs''.

Graphon is not just an intuitive notion, it has mathematical relevance.  This
setting extends work of Aldous~\cite{Aldous1981} and Hoover~\cite{Hoover1979}
in probability theory on exchangeable random graphs (see, e.g.,
\cite{Austin2008}). Graphons are also not just a generalization. They present a
convenient and useful way to study extremal problems for graphs (such as to
find maximum number of edges of a graph satisfying given local properties).

These problems then often take the form of linear inequalities. Lovász
introduced graph algebras (of ``quantum graphs'') with nice ``pictorial''
proofs, see~\cite{Lovasz2012}), and independently Alexander Razborov developed
``flag algebras''~\cite{Razborov2007} which proved to be a very efficient tool
in various extremal problems, see, e.g.,
\cite{HatamiHladkyKralNorinRazborov2012} and~\cite{GrzesikKralLovasz2020}.

The graph algebra of Lovász and Razborov was motivated by early examples
provided by the Caccetta--Häggkvist conjecture, see~\cite{Bondy1997}, the
Sidorenko conjecture~\cite{Sidorenko1991}, and the early
paper~\cite{ErdosLovaszSpencer1979} of Erdős, Lovász, and Spencer on
topological properties of the graphcopy function.

A typical extremal problem may be expressed as a fact that a certain linear
inequality built from homomorphism densities of graphs is nonnegative. This in
turn led Lovász to a question whether any such inequality can be deduced from a
sum of squares of ``quantum graphs''. A related question was formulated by
Razborov~\cite{Razborov2007} whether the validity of any such inequality can be
solved by ``Cauchy--Schwarz Calculus''. However, Hamed Hatami and Serguei
Norin~\cite{HatamiNorin2011} showed that both these questions have a negative
answer in general as the related problems are algorithmically undecidable. So,
extremal problems may be more difficult as originally thought. This was further
supported by the universality results of Cooper, Grzesik, Král, Martins, and
L.\,M.~Lovász, see~\cite{CooperKralMartins2018}
and~\cite{GrzesikKralLovasz2020}, claiming particularly that every graphon may
be extended to a ``finitely forcible'' graphon.

This approach also provides an understanding of the celebrated Szemerédi
regularity lemma. The Szemerédi regularity lemma in this interpretation
means an approximation of every graph (and every graphon) by means of a
``small'' pixel image where almost all entries are constant (but may be
different for different pixels).

The key of the approach of \cite{BorgsChayesLovaszSosVesztergombi2008} is to
characterize convergence using the \emph{cut metric} $d_{\square}(G, H)$ (based
on the cut norm introduced by R.~Frieze and R.~Kannan
in~\cite{FriezeKannan1999}). If the homomorphism density is defined by scaled
subgraph density, then the cut metric is, somewhat dually, characterized by
means of a scaled density of partitions.

The cut metric $d_{\square}(G,H)$ for finite graphs $G,H$ on the
same vertex set $V$ is defined as
\[
\max_{S,T \subseteq V}
\frac{|e_{G}(S,T)| - |e_{H}(S,T)|}
{|V\times V|}
\]
i.e., as the scaled difference of the-sizes of cuts in $G$ and $H$; above
$e_{G}(S,T)$ is the number of edges of $G$ between sets $S$ and $T$. (This
definition can be extended to graphs on different vertex sets. This is
technical and it takes three full pages in~\cite{Lovasz2012}). Interestingly,
the cut distance for a graphon $W$ is more easily defined than in the finite
case: it is induced by the norm:
\[
\|W\| = 
\sup_{S,T \subseteq [0,1]} \,\, \int\limits_{S \times T} W(x,y)dx dy 
\]
The cut norm is also very natural and fitting from an algorithmic point of
view; and it is bounded by the Grothendieck norm up to a multiplicative
constant (as shown by Alon and Naor~\cite{AlonNaor2006}).

As a culmination of several auxiliary results, one obtains that the convergence
is indeed induced by a distance. This is the key fact in many applications and
was proved by Lovász and Szegedy in~\cite{LovaszSzegedy2006}:

\begin{nntheorem}
If $(G_{n})$ is a sequence of graphs of unbounded size, then $(G_{n})$ is a
converging sequence if and only if $(G_{n})$ is a Cauchy sequence with respect
to cut distance $d_{\square} (G_{i},G_{j})$.
\end{nntheorem}

The following result was proved by Lovász and Szegedy
in~\cite{LovaszSzegedy2007}. Lovász considers it as one of the basic results
treated in his book~\cite{Lovasz2012}.

\begin{nntheorem}
The space of all graphons $W$ with cut distance is compact.
\end{nntheorem}

This compactness theorem may be viewed as the roof result for the Szemerédi
regularity Lemma and its various extensions. It also displays the usefulness of
the limit language and of the much more general setting. This area was studied
extensively, for instance by Borgs, Chayes, Elek, Lovász, Sós, Szegedy,
Vesztergombi, and Tao in~\cite{BorgsChayesLovaszSosVesztergombi2012,
ElekSzegedy2012, LovaszSzegedy2007, Tao2006}.

The mathematical richness of this area is best illustrated by the Appendix~A
of~\cite{Lovasz2012} which contains the following sections: Möbius functions;
the Tutte polynomial; some background in probability and measure theory;
moments and the moment problem; ultraproduct and ultralimit;
Vapnik--Chervonenkis dimension; nonnegative polynomials; categories.
Obviously, it is impossible to present here more than a glimpse of what the
book~\cite{Lovasz2012} covers.

Note that the above results are interesting for dense graphs. For sparse graphs
(for example for graphs with constant degrees) one has to devise a different
approach. Limit objects are now called graphings and modelings. For them
results similar to above three theorems are not known. This is treated, e.g.,
by Benjamini and Schramm~\cite{BenjaminiSchramm2001} and by Nešetřil and Ossona
de Mendez~\cite{NesetrilOssona2020}; see again~\cite{Lovasz2012}.

It is amazing that the area of graphs and their limits can be traced back to
Lovász's very early algebraic results (mentioned in Section~\ref{sec:GN2}).
Some forty years later it blossomed in the inspiring climate of the Microsoft
Research Theory Group at Redmond in an atmosphere of concentrated research and
quality, with persons such as Michael Freedman, Oded Schramm and many other
great visitors and with László Lovász as a driving force.

\section{Final Remarks}
\label{sec:GN18}

Let us finish the fireworks of beautiful theorems ranging over many parts of
mathematics and theoretical computer science by adding a few general remarks.

It happens very rarely that a well-known and long-standing open problem is
solved by a novel technique that immediately influences not just that area, but
other parts of mathematics as well. Lovász not only accomplished this once. It
is unbelievable that Lovász repeatedly offered to the world community exactly
such solutions. Some of these proofs are really elegant and were included in
the collections of other beautiful ``book proofs'',
see~\cite{AignerZiegler1998} and~\cite{Matousek2010}.

In this article we concentrated on Lovász-results which had general influence,
led to intensive research by many others, and sometimes spawned the emergence
of whole new theories. Work in areas such as combinatorial optimization,
applications of the ellipsoid method, algebraic graph theory, graph
homomorphisms, topological graph theory, and graph limits is very difficult to
imagine without the pioneering accomplishments of László Lovász.

\begin{figure}[b]
\includegraphics[width=\linewidth]{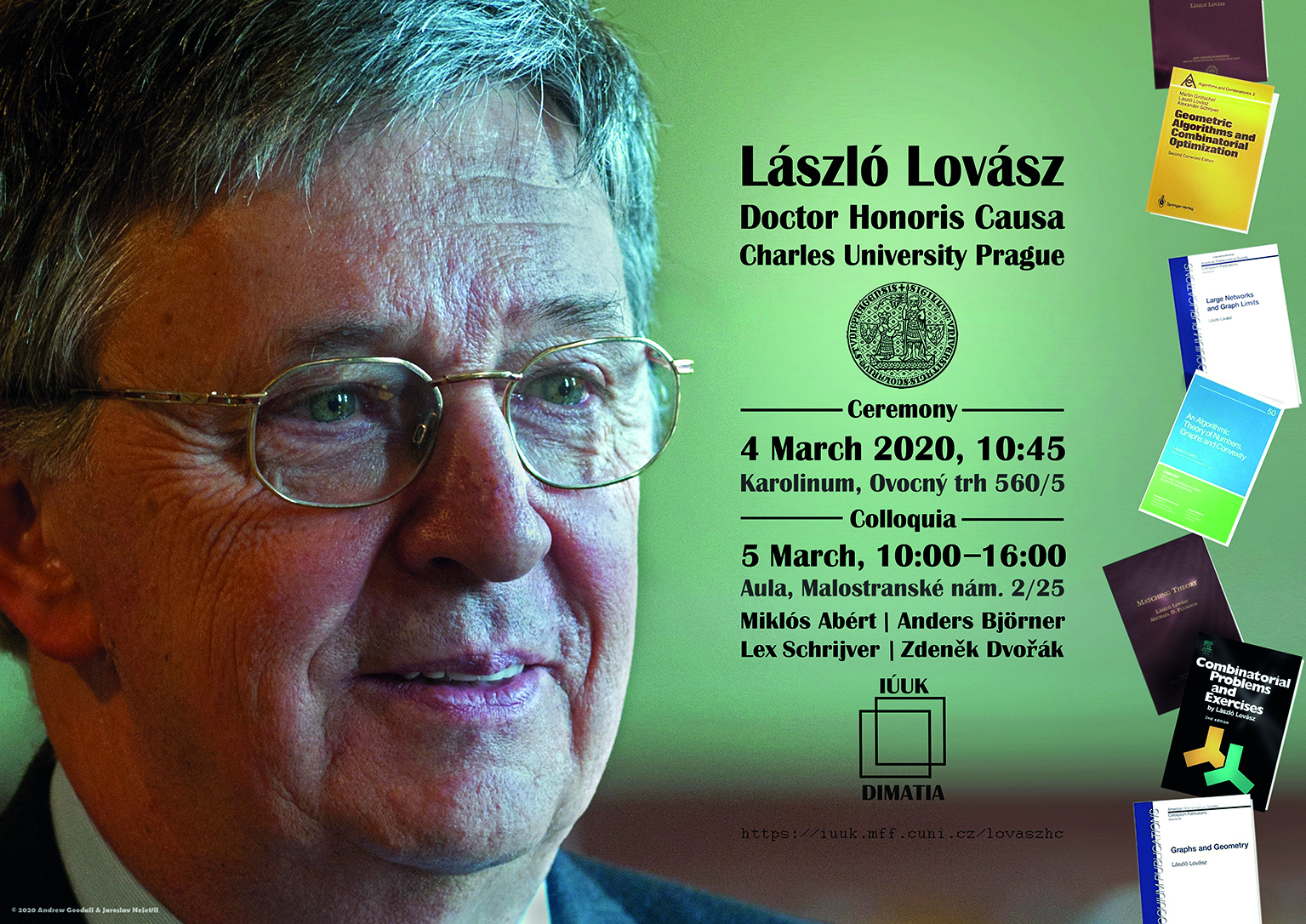}

\caption{Several Lovász-books on a poster (by A.~Goodall and J.\,N.) of the Charles
University in Prague (Photo: Private)}
\label{fig:GN15}

\end{figure}

In our Introduction we indicated that Lovász is both, a ``problem solver'' and
a ``theory builder'', and pointed out that the trio \emph{depth, elegance, and
inspiration} is a particular signature of his work that makes his achievements
unique. We do hope that the glimpse into his oeuvre and the scientific
influence of his results, that we have offered here, provides at least a
partial proof of our conviction.

To keep this article at a reasonable length we had to omit many topics
on which Lovász left his marks. In particular, it was impossible to give
adequate attention to the books he has written, see Fig. 15, and the
influence they had and still have. To mend this omission, albeit very
incompletely, we elucidate the contents and impact of four of his books
-- extremely briefly, though.

Lovász's third book \emph{Combinatorial Problems and
Exercises}~\cite{Lovasz1979b} became \mbox{-- without} any exaggeration -- a
bible for combinatorialists worldwide. This is a book organized in an unusual
way.  It has three parts: The first part consists of mostly easily formulated
questions and problems, the second part contains hints for the solutions, and
the third part thorough proofs with discussions. This of course, makes up the
largest part.

Lovász convincingly claims in this book that discrete mathematics, at the time
of publication, has grown out of an area with simple questions that are
relatively easy to solve without much mathematical knowledge into a structured
field with various branches consisting of central concepts and theorems forming
a hierarchy and possessing a rich bouquet of proof techniques. Instead of
presenting the theories analytically and deductively, Lovász designed his book
with the purpose of helping interested readers to learn many of the existing
techniques in combinatorics. And as he wrote in the introduction: 
\begin{quote}
\emph{The most effective (but admittedly very-time consuming) way of learning
such techniques is to solve (appropriately chosen) exercises and problems}.
\end{quote}
We
believe that this book significantly changed the level on which combinatorics
(and graph theory in particular) was treated. It caught worldwide attention
from the very start (see, e.g., the book review by
Bollobás~\cite{Bollobas1981}) by combinatorialists, computer scientists, and
mathematicians in general. It is remarkable that after more than 40 years of
its existence the book, that mirrors the vast experience of the author, is
still in print and in use.

A side remark: Combinatorics meetings usually have an open problems session
where participants explain questions they are working on and have not solved
yet. Lovász, with his wide knowledge of proof techniques, has always been
outstanding in being able to solve many of the open problems on the spot.

Matching problems have played a considerable role in the development of graph
theory. Well-known and important early results are, e.g., König's Matching
Theorems, the Marriage Theorem, and Tutte's $f$-factor theorem.  Matchings,
$b$-matchings, $T$-joins, etc. have a rich structure theory. The Edmonds--Galai
decomposition is one such example. Various matching problems and their
ramifications appear in a large variety of applications of combinatorial
optimization (e.g., the Chinese Postman Problem). Many of these are solvable
with (highly nontrivial) polynomial time algorithms for which the pioneering
work of J.~Edmonds, see~\cite{Edmonds1965a}, laid the basis.
Edmonds~\cite{Edmonds1965b} achieved also a breakthrough in polyhedral
combinatorics by providing a linear description of the matching polytope that
does not simply follow from total unimodularity. Lovász~\cite{Lovasz1979c} came
up with a new and elegant proof of this result that was later often mimicked
for the characterization of other polytopes arising in combinatorial
optimization.

The book~\cite{LovaszPlummer1986} \emph{Matching Theory}, written by László
Lovász and Mike Plummer, provides a broad view of this subject and covers the
roughly 40~articles that Lovász has contributed to this field. We just want to
highlight Chapters~10 and 11 of this book.  Chapter 10 is devoted to the
$f$-factor problem which asks whether, for a given graph $G=(V,E)$ and integers
$f(v)$ for every vertex $v \in V$, there is a spanning subgraph $H$ of $G$ such
that the degree of $v$ in $H$ is equal to $f(v)$. In a series of four papers
that appeared 1970--1972, Lovász developed a generalization of the
Edmonds--Galai Structure Theorem to the $f$-factor problem to provide an
elegant answer of the $f$-factor problem.  Chapter~11 introduces further
generalizations such as the matroid and polymatroid matching problem which are
interesting (and difficult) combinations of topics in graph and matroid theory.
We refer to this Chapter of~\cite{LovaszPlummer1986} and the
article~\cite{Lovasz1980} for some of the results that can be shown in this
context. Finally, this book contains in the preface a wonderful brief, yet
in-depth survey of the historical development of matching theory.

Lovász's book \emph{Large networks and graph limits}~\cite{Lovasz2012} is
aiming in a different direction. It is the result of a stay of Lovász at the
IAS in Princeton. We have dealt with parts of this book in Section~17. Graph
limits became a very active field with contributions ranging from model theory,
probability, functional analysis to theoretical computer science, network
science and, of course, combinatorics. This theory fits very well with advanced
combinatorics; for example, the role of Szemerédi's regularity lemma is
highlighted and explained properly in this context. The basic theory of
convergent graph sequences is derived in several settings; and multiple
applications to parameter and property testing, extremal theory, and other
applications are given. The book starts with an informal introduction into
large graphs in a network science context, specifying the abundance of real
applications, and questions to ask about them.  This is followed by a lengthy
chapter on the algebra of graph homomorphisms. This chapter can be read
independently and is also of independent interest. But one of the main features
of this book is to show how this algebra is connected to limit structures and
limit distributions. It is amazing how much material was developed in this
context in less than a decade. In the very nice preface, Lovász lists the
branches of mathematics that come into play in his book and writes:

\begin{quote}
\emph{These connections with very different parts of mathematics made
it quite difficult to write this book in a readable form} [$\ldots$] [continuing
that he found that] \emph{the most exciting feature of this theory}
[$\ldots$] [is] \emph{its rich connections with other parts of mathematics
(classical and non-classical)} [$\ldots$] [so that he] \emph{decided to explain
as many of these connections} [$\ldots$] [as he] \emph{could fit in the book}.
\end{quote}
Summarizing, this book is a real tour de force.

The American Mathematical Society Colloquium Publications were established
in~1905. So far 66~books were published in this AMS flagship book series
``\emph{offering the finest in scholarly mathematical publishing}''. Vol.~60 is
the book~\cite{Lovasz2012} \emph{Large Networks and Graph Limits} discussed
above, Vol.~65 is the book \emph{Graphs and Geometry}~\cite{Lovasz2019}, so far
the last book written by Lovász.

Vol. 60 pictures the emergence and maturation of a new theory while Vol.
65 presents a wide spectrum of geometry related techniques (and tricks)
to study graphs. In twenty chapters (and three appendices) Lovász
surveys many connections between graph theory and geometry concentrating
on those which lie deeper. These are among others: rubber band
representations, coin representations, orthogonal representation, and
discrete analytic functions. Interestingly, this book is only about
geometry, and thus topology is outside its scope. Nevertheless, the book
contains some of the key discoveries of Lovász in a new context.

The Leitmotiv of the whole book~\cite{Lovasz2019} is described in the
preface: 

\begin{quote}
\emph{Graphs are usually represented as geometric objects drawn in the plane,
  consisting of vertices and curves connecting them.  The main message of this
  book is that such a representation is not merely a way to visualize the
  graph, but an important mathematical tool.  It is obvious that this geometry
  is crucial in engineering if you want to understand rigidity of frameworks
  and mobility of mechanisms. But even if there is no geometry directly
  connected to the graph-theoretic problem, a well-chosen geometric embedding
has mathematical meaning and applications in proofs and algorithms. This
thought emerged in the 1970s, and I found it quite fruitful}.
\end{quote}
Lovász has been developing these thoughts for about forty years observing: 

\begin{quote}
\emph{Many new results
and new applications of the topic have also been emerging, even outside
mathematics, like in statistical and quantum physics and computer
science (learning theory). At some point I had to decide to round things
up and publish this book}.
\end{quote}
This finishes his preface. But he returns to these considerations in Chapter~20,
``Concluding Thoughts'', on page~390 as follows: 

\begin{quote}
\emph{I am certain that many new results of this nature will be obtained in the
  future (or are already in the literature, sometimes in a quite different
  disguise). Whether these will be collected and combined in another monograph,
  or integrated into science through some other platform provided by the fast
changing technology of communication, I cannot predict. But the beauty of
nontrivial connections between combinatorics, geometry, algebra and physics
will remain here to inspire research}.
\end{quote}
When reviewing the book~\cite{Lovasz2012} in the Bulletin of the American
Mathematical Society, one of us quoted Michel Mendès France who once told him
that envy is the right feeling when reading beautiful mathematics. Yes, this is
the feeling one may have when reading Lovász's books such as~\cite{Lovasz2012}
and~\cite{Lovasz2019}.

His exceptional research capabilities and his broad knowledge of mathematics
are mirrored in Lovász's public presentations and survey articles. He has the
ability to explain difficult results in understandable language and, in
particular, to display and illustrate connections between seemingly unrelated
topics. Examples of that can, e.g., be found in the articles he contributed to
the \emph{Handbook of Combinatorics}~\cite{GrahamGroetschelLovasz1995}, see
also~\cite{Laczkovich2021}. The titles of some of his survey and motivating
articles contain phrases such as \emph{One mathematics} or \emph{Discrete and
Continuous: Two sides of the same}. This reflects his philosophy that science
is not a collection of independent topics but a tightly connected network to be
discovered and understood. He contributed to this conviction also
administratively by serving the scientific community in leading positions of
the International Mathematical Union and the Hungarian Academy of Sciences.

The unity of mathematics and the role of mathematics in the world have been
addressed again and again by László Lovász through many of his activities.
Given the outstanding excellence in his own research and the huge experience as
a professional in combination with admirable modesty the mathematical community
can hardly think of a better representative.

\end{document}